\documentclass[a4paper,11pt,fleqn]{article}
\usepackage{amsmath,amssymb,amsthm,graphicx,subfigure,float,caption,epstopdf,tabularx,color, bm,amsfonts,epic}
\usepackage[top=1in, bottom=1in, left=1.25in, right=1.25in]{geometry}
\usepackage{appendix}
\allowdisplaybreaks
\newtheorem{thm}{Theorem}[section]
\newtheorem{lem}{Lemma}[section]

\newtheorem{rmk}{Remark}[section]
\newtheorem{shm}{Scheme}[section]
\newtheorem*{prf}{Proof}
\numberwithin{equation}{section}

\usepackage{indentfirst}
\graphicspath{{Fig}}

\begin{document}
\title{Arbitrarily high-order energy-preserving schemes for the Camassa-Holm equation}
\author{Chaolong Jiang$^1$,\ Yushun Wang$^2$, \ Yuezheng Gong$^3$\footnote{Correspondence author. Email:
gongyuezheng@nuaa.edu.cn.}\\
{\small $^1$ School of Statistics and Mathematics, }\\
{\small Yunnan University of Finance and Economics, Kunming 650221, P.R. China}\\
{\small $^2$ Jiangsu Provincial Key Laboratory for NSLSCS,}\\
{\small School of Mathematical Sciences,  Nanjing Normal University,}\\
{\small  Nanjing 210023, P.R. China}\\
{\small $^3$ College of Science,}\\
{\small Nanjing University of Aeronautics and Astronautics, Nanjing 210016, P.R. China}\\
}
\date{}
\maketitle

\begin{abstract}
In this paper, we develop a novel class of arbitrarily high-order energy-preserving schemes for the Camassa-Holm equation. With the aid of the invariant energy quadratization approach, the Camassa-Holm equation is first reformulated into an equivalent system, which inherits a quadratic energy. {The new system is then discretized} by the standard Fourier pseudo-spectral method, which can exactly preserve the semi-discrete energy conservation law. Subsequently, { a symplectic Runge-Kutta
method such as the Gauss collocation method is applied} for the resulting semi-discrete system to arrive at an arbitrarily high-order fully discrete scheme.  {We prove that the obtained schemes can conserve the discrete energy conservation law}.  Numerical results are addressed to confirm accuracy and efficiency of the proposed schemes.
  \\[2ex]
\end{abstract}

\section{Introduction}\label{ch:Sec.1}
In this paper, we consider the following Camassa-Holm (CH)  equation \cite{CH93,CHH94}
\begin{align}\label{ch-eq:1.1}
&u_t-u_{xxt}+3uu_x-2u_xu_{xx}-uu_{xxx}=0,\ (x,t)\in \Omega\times(0,T],
  \end{align}
  with periodic boundary condition
  \begin{align*}u(x+L,t)=u(x,t),\ (x,t)\in \Omega\times[0,T],
  \end{align*}
   and initial condition
  \begin{align*}u(x,0)=u_0(x),\ x\in\Omega,
  \end{align*}
   where $\Omega=[a,a+L]$ is a bounded domain.
The CH equation arises as a model for the propagation of
unidirectional shallow water waves, with $u(x,t)$ representing the height of the
fluid free surface above a flat bottom. Under the periodic boundary condition, the CH equation possesses the following conserved quantities
\begin{align}\label{ch-eq:1.2}
\mathcal{M}=\int_{\Omega} udx,\ \mathcal{I}=\int_{\Omega}(u^2+u_x^2)dx, \ \mathcal{H}=-\frac{1}{2}\int_{\Omega}(u^3+uu_x^2)dx,
\end{align}
where $\mathcal{M}$, $\mathcal{I}$ and $\mathcal{H}$  correspond to mass, momentum and Hamiltonian energy of the original problem, respectively.

{The analysis and numerical solution of the CH equation has been widely investigated.} Constantin and Escher \cite{CE98} showed that the solution of the CH equation could develop singularities at a finite time, even for smooth initial data value with compact support. Li and Olver \cite{LO00} established the local well-posedness of the CH equation in the nonhomogeneous Sobolev space $H^s$ with $s>3/2$. Numerical strategies for solving the CH equation include finite difference methods \cite{CKR08a,HR06}, pseudo-spectral methods \cite{KL05,HGL19}, local discontinuous Galerkin method \cite{XS08}, operator splitting methods \cite{CSW16,FL09}, multi-symplectic methods \cite{COR08,CR11,ZST11} and other effective methods (e.g., see Refs. \cite{CL08,FMO10}). However, the mentioned methods cannot exactly preserve the energy of the original system.

It is well known that the energy conservation is an important property of Hamiltonian partial differential equations (PDEs). In Ref. \cite{MY09}, Matsuo et al. presented an energy-conserving Galerkin scheme for the CH equation. Gong and Wang constructed an energy-preserving wavelet collocation scheme for the CH equation \cite{GW16b}. Other energy-preserving schemes can be found in Refs. \cite{ELS19,Matsuo10}. {However, the existing energy-preserving schemes are only second order in
time,} which can't provide satisfactory solutions in long time simulations with a given large time step. {Therefore,
in order to compute long time accurate solutions, the time step has to be refined,
leading to expensive costs.} To remove such obstacle,  the most commonly used approach is to construct higher-order energy-preserving schemes, which make large marching steps practical while preserving the accuracy. To our best knowledge, there has been no reference considering a high-order energy-preserving scheme for the CH equation.

In \cite{QM08}, Quispel and McLaren derived third- and fourth-order averaged vector flied (AVF) methods. Subsequently, Li, Wang and Qin \cite{LWQ14} have extended the fourth-order AVF method to sixth order. At present, high-order AVF methods have been successfully applied to develop high-order energy-preserving schemes for Hamiltonian PDEs (e.g., see \cite{CHWG15,JCWL17}). For linear problems, the proposed high-order schemes are concise and a customized fast solver has been presented to solve the resulting discrete linear equations efficiently \cite{JCWL17}. However, since high-order AVF methods require to calculate high-order derivatives of the vector field, the resulting schemes are tedious for a nonlinear problem such as the CH equation. Based on the method of discrete line integral, Brugnano et al. \cite{BIT10} proposed Hamiltonian Boundary Value Methods (HBVMs), which can be recast as a multistage Runge-Kutta (RK) method. HBVMs are of arbitrarily high order and can exactly preserve energy of polynomial Hamiltonian systems. For non-polynomial cases, a practical energy-preserving scheme is usually gained in the sense that the energy error remains bounded within machine precision, but one has to increase RK stages \cite{BI16}. Recently, energy-preserving continuous stage Runge-Kutta (CSRK) methods have been attracting a lot of interest \cite{H10,LW16,Miyatake14,TS12}. This kind of methods can eliminate the limit of HBVMs to cover non-polynomial Hamiltonian systems. However, the implementation of the CSRK method requires the computation of integrals. If the appearing integrals are replaced by a fixed high-order quadrature, the resulting method is closely related to HBVMs \cite{H10}. In addition, the above mentioned high-order energy-preserving methods are only valid for Hamiltonian systems with constant skew-symmetric structural matrix. For Hamiltonian systems with non-canonical structure matrix, these methods should be further discussed (e.g., see \cite{BCMR12,CH11,WW18}).

In this paper, we present a novel strategy for efficiently developing arbitrarily high-order energy-preserving schemes for the CH equation. We first utilize the invariant energy quadratization (IEQ) approach \cite{GZYW18,YZW17,YZWS17,ZYGW17} to transform the CH equation into an equivalent system, which admits a quadratic energy. Then, the reformulated system is discretized in space by the standard Fourier pseudo-spectral method. {Subsequently, arbitrarily high-order fully discrete schemes are derived by applying in time a symplectic RK method such as the Gauss
collocation method.} The newly proposed schemes have several
desired advantages: {(i) they are energy-preserving for the reformulated Camassa-Holm equation and can reach arbitrarily high-order;} (ii) the required RK stages are optimal; {(iii) they can be directly applied for efficiently solving non-canonical Hamiltonian PDEs where the energies are quadratic.} {It is worth noting that the proposed methods preserve the quadratic energy of the
modified system, but not the Hamiltonian energy of the original CH equation.}

The rest of the paper is organized as follows. In Section \ref{ch:Sec.2}, an equivalent reformulation of the CH equation is presented based on the IEQ approach. In Section \ref{ch:Sec.3}, a semi-discrete system, which inherits a quadratic energy, is obtained by using the standard Fourier pseudo-spectral method for the spatial discretization. In Section \ref{ch:Sec.4}, we apply the Gauss collocation method for the semi-discrete system to arrive at a fully discrete scheme, which is proved to be energy-preserving. Several numerical experiments are presented in Section \ref{ch:Sec.5}. We draw some conclusions in Section \ref{ch:Sec.6}.

\section{Model reformulation using the IEQ approach}\label{ch:Sec.2}
In this section, we reformulate the CH equation into an equivalent form with a quadratic energy functional, which is called the IEQ reformulation. The IEQ reformulation for the CH equation provides an elegant platform for efficiently developing arbitrarily high-order energy-preserving schemes.

Firstly, Eq. \eqref{ch-eq:1.1} can be rewritten equivalently into the following Hamiltonian system
\begin{align}\label{ch-eq:2.1}
\frac{\partial u}{\partial t}=\mathcal{D}\frac{\delta \mathcal{H}}{\delta u},
\end{align}
where $\mathcal{D}=(1-\partial_{xx})^{-1}\partial_x$ is a skew-adjoint operator, $\mathcal{H}$ is the Hamiltonian energy
\begin{align}\label{ch-Hamiltonian-energy}
\mathcal{H} = \int_{\Omega}H(u,u_x)dx, \ H(u,u_x) = -\frac{1}{2}(u^3+uu_x^2),
\end{align}
and $\frac{\delta \mathcal{H}}{\delta u}$ denotes the variational derivative of $\mathcal{H}$ with respect to $u$
\begin{align*}
\frac{\delta \mathcal{H}}{\delta u}=\frac{\partial H }{\partial u}-\frac{\partial}{\partial x}\frac{\partial H }{\partial u_x} = -\frac{3}{2}u^2-\frac{1}{2}u_x^2+(u u_x)_x.
\end{align*}
One intrinsic property of \eqref{ch-eq:2.1} is energy-preserving, i.e.,
\begin{align}\label{hg-ep:Eq:2.2}
\frac{d}{dt}\mathcal{H}=(\frac{\delta \mathcal{H}}{\delta u},\frac{\partial u}{\partial t})=(\frac{\delta \mathcal{H}}{\delta u},\mathcal{D}\frac{\delta \mathcal{H}}{\delta u})=0,
\end{align}
where $(\cdot,\cdot)$ is the $L^2$-inner product defined by $(f,g)=\int_{\Omega} fg dx$.

Then, we formulate the idea of the IEQ approach for the Hamiltonian system \eqref{ch-eq:2.1}. Let
\begin{align*}
q=g(u,u_x)=-\frac{1}{2}(u^2+u_x^2),
\end{align*}
the Hamiltonian energy functional is then rewritten as
\begin{align}\label{ch-eq:2.3}
\mathcal{H}=\int_{\Omega}uqdx.
\end{align}
According to energy variational, the Hamiltonian system  \eqref{ch-eq:2.1} can be reformulated into the following equivalent form
\begin{align}\label{ch-eq:2.4}
\left\lbrace
  \begin{aligned}
  &\partial_t u=\mathcal{D}\Bigg(q+u\frac{\partial g}{\partial u}-\partial_x\big(u\frac{\partial g}{\partial u_x}\big)\Bigg), \\
  &\partial_tq=\frac{\partial g}{\partial u}u_t+\frac{\partial g}{\partial u_x}\partial_xu_{t},
  \end{aligned}\right.\ \ 
  \end{align}
with consistent initial conditions
  \begin{align*}
  u(x,0)=u_0(x),\ q(x,0)=-\frac{1}{2}\big(u(x,0)^2+u_x(x,0)^2\big),
  \end{align*}
where
  \begin{align*}
  \frac{\partial g}{\partial u}=-u,\ \frac{\partial g}{\partial u_x}=-u_x.
  \end{align*}

\begin{thm} \label{ch-thm-2-1} The system \eqref{ch-eq:2.4} satisfies the following quadratic energy
\begin{align*}
\frac{d}{dt}\mathcal{H}=0,\ \mathcal{H}=\int_{\Omega}uqdx.
\end{align*}
\end{thm}
\begin{prf}\rm
By some calculations, we obtain from the system \eqref{ch-eq:2.4}
\begin{align*}
\frac{d}{dt}\mathcal{H}&=(u_t,q)+(u,q_t)\\
&=(u_t,q) + \Big(u, \frac{\partial g}{\partial u} u_t + \frac{\partial g}{\partial u_x} \partial_x u_{t}\Big)\\
&=(q,u_t) + (u\frac{\partial g}{\partial u}, u_t) - \Big(\partial_x(u\frac{\partial g}{\partial u_x}), u_{t}\Big)\\
&=\Big(q+u\frac{\partial g}{\partial u}-\partial_x\big(u\frac{\partial g}{\partial u_x}\big), u_t\Big)\\
&=\Bigg(q+u\frac{\partial g}{\partial u}-\partial_x\big(u\frac{\partial g}{\partial u_x}\big), \mathcal{D}\Big(q+u\frac{\partial g}{\partial u}-\partial_x\big(u\frac{\partial g}{\partial u_x}\big)\Big)\Bigg)\\
&=0,
\end{align*}
where the last equality follows from the skew-adjoint of $\mathcal{D}$. This completes the proof.\qed
\end{prf}

\section{Energy-preserving spatial discretization}\label{ch:Sec.3}
Many energy-preserving schemes have been designed and investigated for solving the
CH equation in the literature, but little attention is paid to the energy-preserving
properties brought by spatial discretization. In this section, the Fourier pseudo-spectral method is applied for the reformulated Camassa-Holm equation \eqref{ch-eq:2.4} to derive a spatial semi-discrete scheme, which is shown to preserve the semi-discrete quadratic energy \eqref{ch-eq:2.3}.

Let $\Omega_{h}=\{x_{j}|x_{j}=a+jh,\ 0\leq j\leq N\}$ be a partition of $\Omega = [a,a+L]$ with mesh size $h=L/N$, where $N$ is an even number. A discrete mesh function ${U}_{j} = U(x_j),\ j\in \mathbb{Z}$ satisfies the periodic boundary condition if and only if
\begin{align}\label{ch-PBS}
U_{j}=U_{j+N}.
\end{align}
Let $\mathbb{V}_{h}=\big\{{\bm U}|{\bm U}=(U_{0},U_{1},\cdots,U_{N-1})^{T}\big\}$ be the space of mesh functions on $\Omega_h$ {that satisfy the periodic boundary condition \eqref{ch-PBS}}.
We define the discrete inner product as follows
\begin{align*}
&\langle {\bm U},{\bm V}\rangle_{h}=h\sum_{j=0}^{N-1}U_{j}V_{j},\ \forall\  {\bm U},{\bm V}\in\mathbb{ V}_{h}.
\end{align*}
The discrete $L^{\infty}$-norm of ${\bm U}\in\mathbb{V}_h$ is defined as
\begin{align*}
\|{\bm U}\|_{h,\infty}=\max\limits_{0\le j\le N-1}|U_{j}|.
\end{align*}
In addition, we denote $`\cdot$' as the componentwise product of vectors ${\bm U},{\bm V}\in\mathbb{V}_h$, that is,
\begin{align*}
{\bm U}\cdot {\bm V}=&\big(U_{0}V_{0},U_{1}V_{1},\cdots,U_{N-1}V_{N-1}\big)^{T}.
\end{align*}
For brevity, we denote $\underbrace{{\bm U}\cdot...\cdot {\bm U}}_{p}$ as ${\bm U}^p$.

Let $S_{N}=\text{span}\{g_{j}(x),\ 0\leq j\leq N-1\}$ be the interpolation space, {where $g_{j}(x)$ is the trigonometric polynomial of degree $N/2$ }given by
\begin{align*}
  &g_{j}(x)=\frac{1}{N}\sum_{l=-N/2}^{N/2}\frac{1}{a_{l}}e^{\text{i}l\mu (x-x_{j})},
\end{align*}
with $a_{l}=\left \{
 \aligned
 &1,\ |l|<\frac{N}{2},\\
 &2,\ |l|=\frac{N}{2},
 \endaligned
 \right.$ and $\mu=\frac{2\pi}{L}$.
We define the interpolation operator $I_{N}: C(\Omega)\to S_{N}$ as follows \cite{CQ01}
\begin{align*}
I_{N}U(x)=\sum_{j=0}^{N-1}U_{j}g_{j}(x),
\end{align*}
where $U_{j}=U(x_{j})$. Taking the derivative with respect to $x$, and then evaluating the resulting expression at the collocation point $x_{j}$, we have
\begin{align*}
\frac{\partial^{s} I_{N}U(x_{j})}{\partial x^{s}}
&=\sum_{k=0}^{N-1}U_{k}\frac{d^{s}g_{k}(x_{j})}{dx^{s}}=[{\bm D}_{s}{\bm U}]_{j},\ {\bm U}\in\mathbb{V}_h,
\end{align*}
where $j=0,\cdots,N-1$ and ${\bm D}_{s}$ is an $N\times N$ matrix with elements given by
\begin{align*}
({\bm D}_{s})_{j,k}=\frac{d^{s}g_{k}(x_{j})}{dx^{s}},\ j,k=0,1,\cdots,N-1.
\end{align*}
In particular, for first and second derivatives, we have, respectively
\begin{align*}
\frac{\partial I_{N}U(x_{j})}{\partial x}
=[{\bm D}_{1}{\bm U}]_{j},\
\frac{\partial^{2} I_{N}U(x_{j})}{\partial x^{2}}
=[{\bm D}_{2}{\bm U}]_{j},\ j=0,1,\cdots,N-1,
\end{align*}
where ${\bm D}_{1}$ is a real skew-symmetric matrix, and ${\bm D}_{2}$ is a real symmetric matrix. We note that \cite{ST06}
\begin{align}
&{\bm D}_{1}={\bm F}_{N}^{H}\Lambda_1{\bm F}_{N},\ \Lambda_1=\text{i}\mu\text{diag}\big(0,1,\cdots,\frac{N}{2}-1,0,-\frac{N}{2}+1,\cdots,-1\big),\\
&{\bm D}_{2}={\bm F}_{N}^{H}\Lambda_2{\bm F}_{N},\ \Lambda_2=\big[\text{i}\mu\text{diag}(0,1,\cdots,\frac{N}{2}-1,\frac{N}{2},-\frac{N}{2}+1,\cdots,-1)\big]^2,
\end{align}
 where  ${\bm F}_{N}$ is the discrete Fourier transform matrix with elements
$\big({\bm F}_{N}\big)_{j,k}=\frac{1}{\sqrt{N}}e^{-\text{\rm i}jk\frac{2\pi}{N}},$ ${\bm F}_{N}^{H}$ is the conjugate transpose matrix of ${\bm F}_{N}$.

Applying the standard Fourier pseudo-spectral method to the system \eqref{ch-eq:2.4} in space, we have
\begin{align}\label{ch-eq:3.2}
\left\lbrace
  \begin{aligned}
  &\frac{d}{dt}{\bm U}={\bm D}\Bigg({\bm Q} - {\bm U}^2 + {\bm D}_1 \Big(({\bm D}_1 {\bm U}\big)\cdot{\bm U}\Big)\Bigg), \\
  &\frac{d}{dt}{\bm Q}=-{\bm U}\cdot\frac{d}{dt}{\bm U}-\big({\bm D}_1 {\bm U}\big)\cdot\Big({\bm D}_1\frac{d}{dt}{\bm U}\Big),
  \end{aligned}\right.\
\end{align}
where ${\bm D} = ({\bm I} - {\bm D}_2)^{-1}{\bm D}_1$, and ${\bm U},{\bm Q}\in \mathbb{V}_h.$

\begin{thm}
 The system \eqref{ch-eq:3.2} preserves the following semi-discrete quadratic energy
 \begin{align*}
 \frac{d}{dt}E=0,\ E=\Big\langle{{\bm U},\bm Q}\Big\rangle_h.
 \end{align*}
  \end{thm}
 \begin{prf}\rm The proof strictly follows that done for Theorem \ref{ch-thm-2-1}, thus, for brevity, we omit it.
 \end{prf}

\begin{rmk}\label{rmk3.1} If the Fourier pseudo-spectral method is employed for the original system \eqref{ch-eq:2.1}, we can obtain a new semi-discrete scheme
\begin{align}
\frac{d}{dt}{\bm U}={\bm D}\Bigg(-\frac{3}{2} {\bm U}^2 - \frac{1}{2}({\bm D}_1 {\bm U})^2 + {\bm D}_1 \Big(({\bm D}_1 {\bm U}\big)\cdot{\bm U}\Big)\Bigg),
\end{align}
which can also be proved to preserve a semi-discrete Hamiltonian energy
\begin{align*}
\frac{d}{dt} H = 0,\ H=-\frac{h}{2}\sum_{j=0}^{N-1}\Big(U_j^3+U_j\cdot({\bm D}_1{\bf U})_j^2\Big).
\end{align*}
We note that the quadratic energy \eqref{ch-eq:2.3} is only equivalent to the Hamiltonian energy \eqref{ch-Hamiltonian-energy} in continuous sense, but not for the semi-discrete sense.
\end{rmk}

\section{Energy-preserving fully discretized schemes}\label{ch:Sec.4}
 In this section, we first derive a class of high-order energy-preserving schemes by using the Gauss collocation method in time for the IEQ reformulation \eqref{ch-eq:3.2}. Then, we show that, together with the other time integrators \cite{GWW18,JCW18}, the IEQ reformulation \eqref{ch-eq:3.2} also provides an elegant platform for efficiently developing linear-implicitly energy-preserving schemes.
\subsection{High-order energy-preserving schemes}

%

Applying an $s$-stage collocation method to the system \eqref{ch-eq:3.2} in time, we obtain the following scheme.
 \begin{shm}\label{shm4.2} Let $c_1,\cdots, c_s$ be distinct real numbers (usually $0\le c_i\le 1$). For given ${\bm U}^n,{\bm Q}^n\in\mathbb{V}_h$, ${\bm u}(t)$ and ${\bm v}(t)$ are two $N$ dimensional vector polynomials of degree $s$ satisfying, respectively,
 \begin{align}
 &{\bm u}(t_n)={\bm U}^n,\ {\bm v}(t_n)={\bm Q}^n,\\
 &\dot{\bm u}(t_n^i)={\bm D}\Bigg({\bm v}(t_n^i)-{\bm u}(t_n^i)^2 + {\bm D}_1\Big(\big({\bm D}_1{\bm u}(t_n^i)\big)\cdot{\bm u}(t_n^i)\Big)\Bigg),\\
 &\dot{\bm v}(t_n^i)=-{\bm u}(t_n^i)\cdot\dot{\bm u}(t_n^i)-\big({\bm D}_1{\bm u}(t_n^i)\big)\cdot\big({\bm D}_1\dot{\bm u}(t_n^i)\big),
 \end{align}
 where $t_n^i=t_n+c_i\tau, i=1,\cdots,s$. And the numerical solution is defined by ${\bm U}^{n+1}={\bm u}(t_n+\tau)$ and ${\bm Q}^{n+1}={\bm v}(t_n+\tau)$, respectively.
\end{shm}

\noindent Theorem 1.4 on page 31 of Ref. \cite{ELW06} indicates that the collocation method is equivalent to a RK method. If we take $c_1,\cdots, c_s$ as the zeros of the $s$th shifted Legendre polynomial
 \begin{align*}
\frac{d^s}{dx^s}\Big(x^s(x-1)^s\Big),
\end{align*}
Scheme \ref{shm4.2} is called Gauss collocation method and has order $2s$, and the zeros are called Gauss collocation points. Collocation points for Gauss collocation methods of order 4 and 6 are given explicitly in Ref. \cite{ELW06}.

\begin{rmk}{ As pointed out above, the proposed scheme can not preserve the
Hamiltonian energy \eqref{ch-Hamiltonian-energy} of the original Camassa-Holm equation,  but only the quadratic energy of the modified system \eqref{ch-eq:2.4}.}
\end{rmk}

{It is well known that symplectic Runge-Kutta schemes
preserve all the quadratic invariants of the ODE problem (see Refs. \cite{Cooperima87,Sanz-Sernabit88,SCbook94}). Since Gauss collocation schemes of any order are symplectic (see Ref. \cite{Sanz-Sernabit88} and references therein), the following theorem is straightforward.}

 \begin{thm}\label{thm4.1} The $s$-stage Gauss collocation Scheme \ref{shm4.2} is energy-preserving,
i.e., it satisfies the following quadratic energy
\begin{align}\label{ch-eq:4.7}
E^{n+1}=E^{n},\ E^n=\langle{\bm U}^n,{\bm Q}^n\rangle_h, \ n=0,1,\cdots,M-1.
\end{align}
 \end{thm}

\begin{rmk} {Since any other symplectic Runge-Kutta
method would preserve the quadratic energy, other arbitrarily high-order schemes which preserve the discrete quadratic energy \eqref{ch-eq:4.7} can be easily obtained.}

\end{rmk}

\subsection{Linearly-implicit energy-preserving schemes}
In this subsection, a novel, linearly-implicit and energy-preserving scheme
for the CH equation  is obtained by utilizing in time the linearized Crank-Nicolson method for the semi-discrete system \eqref{ch-eq:3.2}. The resulting scheme is denoted by IEQ-LCNS.
 \begin{shm}\label{schmeA1} Applying the linearized Crank-Nicolson method to discretize the semi-discrete system \eqref{ch-eq:3.2} in time, we obtain a fully discretized scheme, as follows:
\begin{align}\label{ch-A:1}
\left\lbrace
  \begin{aligned}
  &\delta_t^+{\bm U}^n={\bm D}\Bigg({\bm Q}^{n+\frac{1}{2}}+\text{diag}\big(-\hat{\bm U}^{n+\frac{1}{2}}\big){\bm U}^{n+\frac{1}{2}}\\
&~~~~~~~~~~~~~~~~~~~~~~~~~~~~~~~~~~~-{\bm D}_1\Big(\text{diag}\big(-{\bm D}_1\hat{\bm U}^{n+\frac{1}{2}}\big){\bm U}^{n+\frac{1}{2}}\Big)\Bigg), \\
  &\delta_t^+{\bm Q}^n=\text{diag}\big(-\hat{\bm U}^{n+\frac{1}{2}}\big)\delta_t^+{\bm U}^n+\text{diag}\big(-{\bm D}_1\hat{\bm U}^{n+\frac{1}{2}}\big){\bm D}_1\delta_t^+{\bm U}^n.
  \end{aligned}\right.\ \ 
  \end{align}
   where ${\delta_t^{+}{\bm U}^n}=\frac{{\bm U}^{n+1}-{\bm U}^n}{\tau},\hat{\bm U}^{n+\frac{1}{2}}=\frac{3{\bm U}^{n}-{\bm U}^{n-1}}{2}$ and ${\bm Q}^{n+\frac{1}{2}}=\frac{{\bm Q}^{n+1}+{\bm Q}^{n}}{2}$ and ${\bm U}^1,{\bm Q}^1$ is the solution of the following equation
   \begin{align}\label{ch-A:2}
\left\lbrace
  \begin{aligned}
  &\delta_t^+{\bm U}^0={\bm D}\Bigg({\bm Q}^{\frac{1}{2}}+\text{diag}\big(-{\bm U}^{0}\big){\bm U}^{\frac{1}{2}}-{\bm D}_1\Big(\text{diag}\big(-{\bm D}_1{\bm U}^{0}\big){\bm U}^{\frac{1}{2}}\Big)\Bigg), \\
  &\delta_t^+{\bm Q}^0=\text{diag}\big(-{\bm U}^{0}\big)\delta_t^+{\bm U}^0+\text{diag}\big(-{\bm D}_1{\bm U}^{0}\big){\bm D}_1\delta_t^+{\bm U}^0.
  \end{aligned}\right.\ \ 
  \end{align}
\end{shm}
\begin{thm}  The {IEQ-LCNS} \eqref{ch-A:1}-\eqref{ch-A:2} can exactly preserve the discrete quadratic energy \eqref{ch-eq:4.7}.

\end{thm}
\begin{prf} \rm We can deduce from \eqref{ch-A:1} that
\begin{align}\label{ch:A3}
\delta_t^+\langle {\bm Q}^n,{\bm U}^n\rangle_h&=\langle \delta_t^+{\bm Q}^n,{\bm U}^{n+\frac{1}{2}}\rangle_h+\langle {\bm Q}^{n+\frac{1}{2}},\delta_t^+{\bm U}^n\rangle_h\nonumber\\
&=\langle\text{diag}\big(-\hat{\bm U}^{n+\frac{1}{2}}\big)\delta_t^+{\bm U}^n+\text{diag}\big(-{\bm D}_1\hat{\bm U}^{n+\frac{1}{2}}\big){\bm D}_1\delta_t^+{\bm U}^n,{\bm U}^{n+\frac{1}{2}}\rangle_h\nonumber\\
&~~~~~+\langle {\bm Q}^{n+\frac{1}{2}},\delta_t^+{\bm U}^n\rangle_h\nonumber\\
&=\langle\text{diag}\big(-\hat{\bm U}^{n+\frac{1}{2}}\big){\bm U}^{n+\frac{1}{2}}-{\bm D}_1\big(\text{diag}\big(-{\bm D}_1\hat{\bm U}^{n+\frac{1}{2}}\big){\bm U}^{n+\frac{1}{2}}\big),\delta_t^+{\bm U}^n\rangle_h\nonumber\\
&~~~~~+\langle {\bm Q}^{n+\frac{1}{2}},\delta_t^+{\bm U}^n\rangle_h\nonumber\nonumber\\
&=\langle {\bm G}^n,\delta_t^+{\bm U}^n\rangle_h=\langle {\bm G}^n,{\bm D}{\bm G}^n\rangle_h=0,
\end{align}
where
\begin{align*}
{\bm G}^n={\bm Q}^{n+\frac{1}{2}}+\text{diag}\big(-\hat{\bm U}^{n+\frac{1}{2}}\big){\bm U}^{n+\frac{1}{2}}-{\bm D}_1\big(\text{diag}\big(-{\bm D}_1\hat{\bm U}^{n+\frac{1}{2}}\big){\bm U}^{n+\frac{1}{2}}\big),
\end{align*}
and the last equality follows from the skew-symmetry of ${\bm D}$.
Thus, according to \eqref{ch:A3}, we have
\begin{align*}
E^{n+1}=E^{n},\ E^n=\langle {\bm Q}^n,{\bm U}^n\rangle_h,\ n=1,2,\cdots,M-1.
\end{align*}
An argument similar to \eqref{ch-A:2} used in \eqref{ch:A3} shows that
\begin{align*}
E^{1}=E^{0}.
\end{align*}
This completes the proof.\qed

\end{prf}

Subsequently, we show that the above scheme can be solved efficiently. Indeed, we first rewrite \eqref{ch-A:1} as
\begin{align}\label{ch-A:8}
\left\lbrace
  \begin{aligned}
  &{\bm U}^{n+\frac{1}{2}}={\bm U}^n+\frac{\tau}{2}{\bm D}\Bigg[{\bm Q}^{n+\frac{1}{2}}+{\bm g}_1^n{\bm U}^{n+\frac{1}{2}}\Bigg], \\
  &{\bm Q}^{n+\frac{1}{2}}={\bm g}_2^n{\bm U}^{n+\frac{1}{2}}+{\bm Q}^n-{\bm g}_2^n{\bm U}^n,
  \end{aligned}\right.\ \ 
  \end{align}
  where
  \begin{align*}
&{\bm g}_1^n=\text{diag}\big(-\hat{\bm U}^{n+\frac{1}{2}}\big)-{\bm D}_1\Big(\text{diag}\big(-{\bm D}_1\hat{\bm U}^{n+\frac{1}{2}}\big)\Big),\\
&{\bm g}_2^n=\text{diag}\big(-\hat{\bm U}^{n+\frac{1}{2}}\big)+\text{diag}\big(-{\bm D}_1\hat{\bm U}^{n+\frac{1}{2}}\big){\bm D}_1.
\end{align*}
Then, by eliminating ${\bm Q}^{n+\frac{1}{2}}$ from \eqref{ch-A:8}, we have
 \begin{align}\label{ch:A6}
 {\bm U}^{n+\frac{1}{2}}=\frac{\tau}{2}{\bm D}\Bigg[{\bm g}_1^n{\bm U}^{n+\frac{1}{2}}+{\bm g}_2^n{\bm U}^{n+\frac{1}{2}}\Bigg]+{\bm b}^n,
 \end{align}
 where
 \begin{align*}
 {\bm b}^n={\bm U}^n+\frac{\tau}{2}{\bm D}\Big({\bm Q}^n-{\bm g}_2^n{\bm U}^{n}\Big).
 \end{align*}
 Finally, we obtain ${\bm U}^{n+\frac{1}{2}}$ from \eqref{ch:A6} by using the following iteration method for linear equations \eqref{ch:A6}
\begin{align}
{\bm U}^{n+\frac{1}{2},s+1}=\frac{\tau}{2}{\bm D}\Bigg[{\bm g}_1^n{\bm U}^{n+\frac{1}{2},s}+{\bm g}_2^n{\bm U}^{n+\frac{1}{2},s}\Bigg]+{\bm b}^n,
\end{align}
where we take the initial iteration vector ${\bm U}^{n+\frac{1}{2},0}={\bm U}^n$ and each iteration will terminate if
the infinity norm of the error between two adjacent iterative steps is less than $10^{-14}$. Then, ${\bm Q}^{n+\frac{1}{2}}$ {is obtained from the second equality of} \eqref{ch-A:8}. Subsequently, we have ${\bm U}^{n+1}=2{\bm U}^{n+\frac{1}{2}}-{\bm U}^n$ and ${\bm Q}^{n+1}=2{\bm Q}^{n+\frac{1}{2}}-{\bm Q}^n$.
\begin{rmk}
We should note  from \eqref{ch:A6} that, the IEQ approach need introduce an auxiliary variable, but the auxiliary variable can be eliminated in practical computations.

\end{rmk}


\section{Numerical examples}\label{ch:Sec.5}
In this section, we will investigate the accuracy, {CPU time and invariants-preservation of the proposed schemes}. Theoretically, the newly proposed high-order schemes \ref{shm4.2} could reach arbitrarily high-order accuracy in time (with proper choice of the Gauss collocation points), and they all can exactly preserve the discrete quadratic energy \eqref{ch-eq:4.7}. For simplicity, in the rest of this paper, the Gauss methods of order 4 (denoted by 4th-order HIEQ-GM) and 6 (denoted by 6th-order HIEQ-GM) are only used for demonstration purposes. Also, the results are compared with the energy-preserving Fourier pseudo-spectral scheme (denoted by EPFPS) and the multi-symplectic Fourier pseudo-spectral scheme (denoted by MSFPS), where we substitute the Fourier pseudo-spectral method into the wavelet collocation method for space directions in Refs. \cite{GW16b} and \cite{ZST11}, respectively. For the convergence rate, we use the formula
\begin{align*}
\text{Rate}=\frac{\ln( error_{1}/error_{2})}{\ln (\tau_{1}/\tau_{2})},
\end{align*}
where $\tau_{l}, error_{l}, (l=1,2)$ are step sizes and errors with the step size $\tau_{l}$, respectively.
\subsection{Accuracy test}
\begin{figure}[H]
\centering\begin{minipage}[t]{80mm}
\includegraphics[width=85mm]{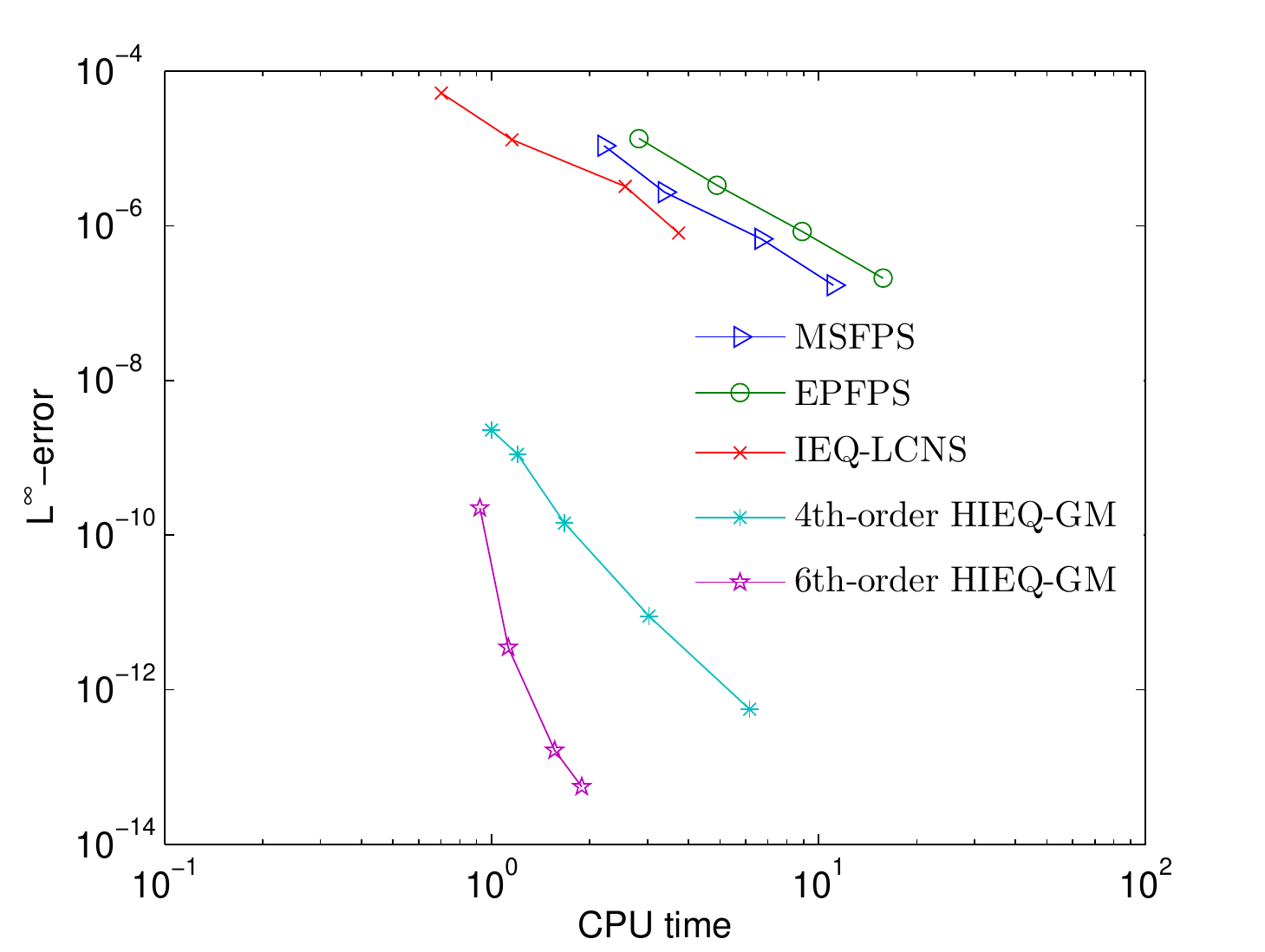}
\end{minipage}
\caption{The numerical error versus the CPU time.}\label{Fig_ch:1}
\end{figure}

We consider the periodic smooth solution with initial condition
\begin{align*}
u_0(x)=\sin(x),\ x\in\Omega=[0,2\pi],
\end{align*}
and the periodic boundary condition. The \emph{exact} solution is obtained numerically by 6th-order HIEQ-GM under a
very small time step $\tau=0.001$ and spatial step $h=\frac{2\pi}{128}$ at $T=1$.  In Tables \ref{Tab-ch:1} and \ref{Tab-ch:2}, we display the numerical error in discrete $L^{\infty}$-norm and the convergence rate for different schemes at $T=1$, respectively.
As illustrated in Table \ref{Tab-ch:1}, all of the schemes have second-order convergence rate in time, and from {Table \ref{Tab-ch:2}}, it is clear to see that 4th-order HIEQ-GM and 6th-order HIEQ-GM can arrive at fourth-order and sixth-order convergence rates in time, respectively. In Fig. \ref{Fig_ch:1}, we plot the global numerical error in discrete $L^{\infty}$-norm versus the CPU time for different schemes. The plot shows
that, for a given global error, the sixth-order scheme is computationally cheapest. The IEQ-LCNS admits larger numerical errors than the ones provided by EPFPS and MSFPS, however, it is computationally cheaper.
\begin{table}[H]
\tabcolsep=9pt
\footnotesize
\renewcommand\arraystretch{1.1}
\centering
\caption{{The numerical error and convergence rate for different second-order schemes with $h=\frac{2\pi}{128}$ and different time steps at $T=1$.}}\label{Tab-ch:1}
\begin{tabularx}{1.0\textwidth}{XXXXXXXXX}\hline
{Scheme\ \ } &{$\tau$} &{$L^{\infty}$} &{Rate}\\
\hline
 {IEQ-LCNS}&{$\frac{1}{100}$} & {2.083e-04}&{-}\\[1ex]
 {}&{$\frac{1}{200}$}& {5.182e-05}&{ 2.01}\\[1ex]
  {}  &{$\frac{1}{400}$}& {1.293e-05}&{2.00}\\[1ex]
  {}  &{$\frac{1}{800}$}& {3.230e-06}&{2.00}\\[1ex]
    {EPFPS} &{$\frac{1}{100}$} & { 4.355e-05}&{-}\\[1ex]
   {} &{$\frac{1}{200}$} & { 1.089e-05}&{2.00}\\[1ex]
 {}  &{$\frac{1}{400}$}& {2.722e-06}&{2.00} \\[1ex]
  {}  &{$\frac{1}{800}$}& {6.806e-07}&{2.00} \\[1ex]
  {MSFPS} &{$\frac{1}{100}$} & {4.344e-05}&{-}\\[1ex]
   {} &{$\frac{1}{200}$} & {1.086e-05}&{2.00}\\[1ex]
 {}  &{$\frac{1}{400}$}& { 2.716e-06}&{2.00} \\[1ex]
  {}  &{$\frac{1}{800}$}& { 6.790e-07}&{2.00} \\\hline
\end{tabularx}
\end{table}

\begin{table}[H]
\tabcolsep=9pt
\footnotesize
\renewcommand\arraystretch{1.1}
\centering
\caption{{The numerical error and convergence rate for different high-order schemes with $h=\frac{2\pi}{128}$ and different time steps at $T=1$.}}\label{Tab-ch:2}
\begin{tabularx}{1.0\textwidth}{XXXXXXXXX}\hline
{Scheme\ \ } &{$\tau$} &{$L^{\infty}$} &{Rate}\\ 
\hline
 {4th-order HIEQ-GM} &{$\frac{1}{30}$} & {2.817e-07}&{-} \\ [1ex]
 {}  &{$\frac{1}{60}$}& {1.765e-08}&{4.00}  \\[1ex]
  {}  &{$\frac{1}{120}$}& {1.104e-09}&{4.00}  \\[1ex]
 {6th-order HIEQ-GM} &{$\frac{1}{30}$} & {2.231e-010}&{-}\\[1ex]
 {} &{$\frac{1}{60}$}& {3.523e-012}&{5.98}  \\[1ex]
  {}  &{$\frac{1}{120}$}& { 5.534e-014}&{5.99}  \\\hline
\end{tabularx}
\end{table}

\subsection{Peakon solution}
\begin{figure}[H]
\centering\begin{minipage}[t]{60mm}
\includegraphics[width=65mm]{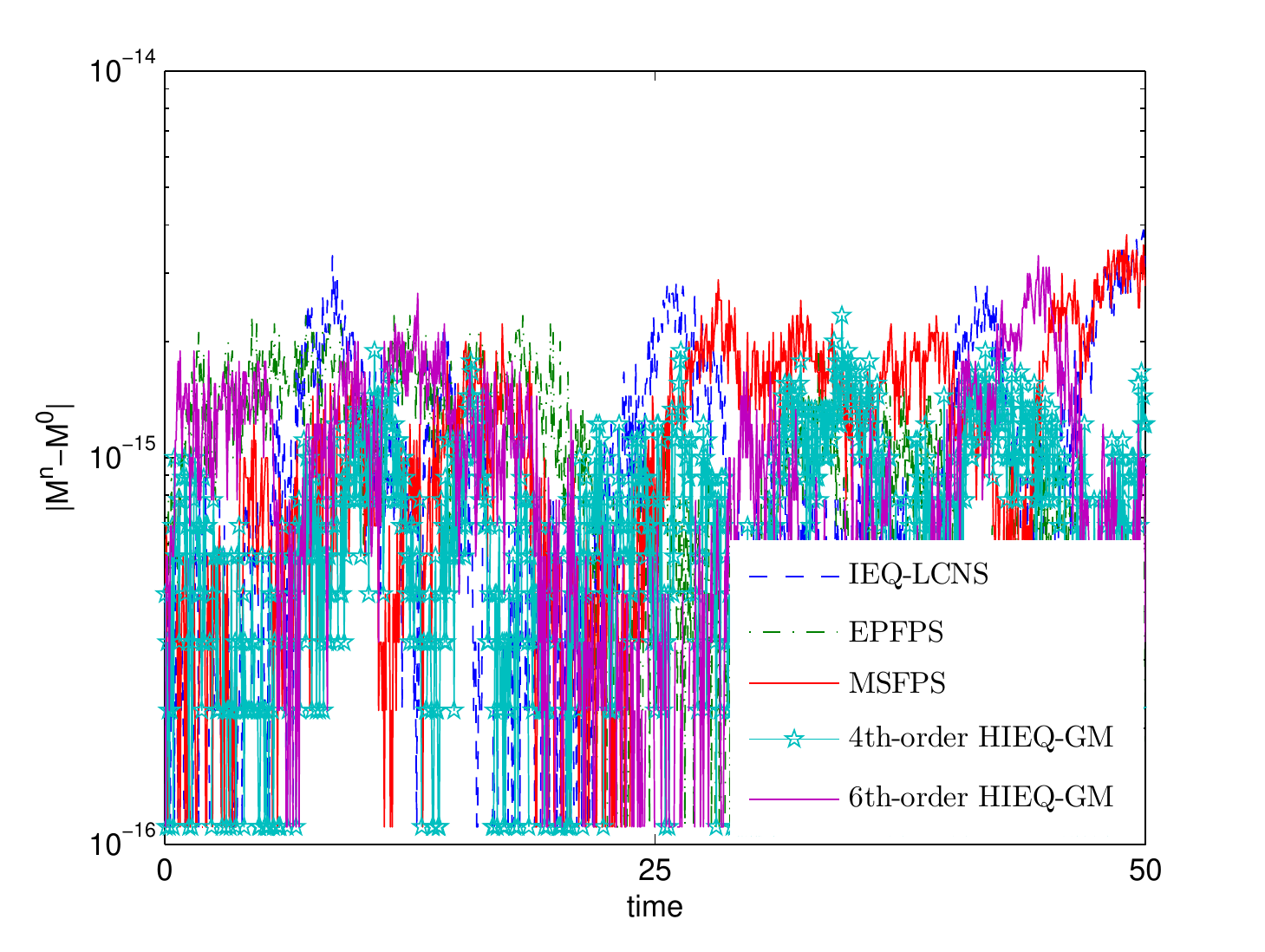}
\caption*{(c) Mass}
\end{minipage}\ \
\begin{minipage}[t]{60mm}
\includegraphics[width=65mm]{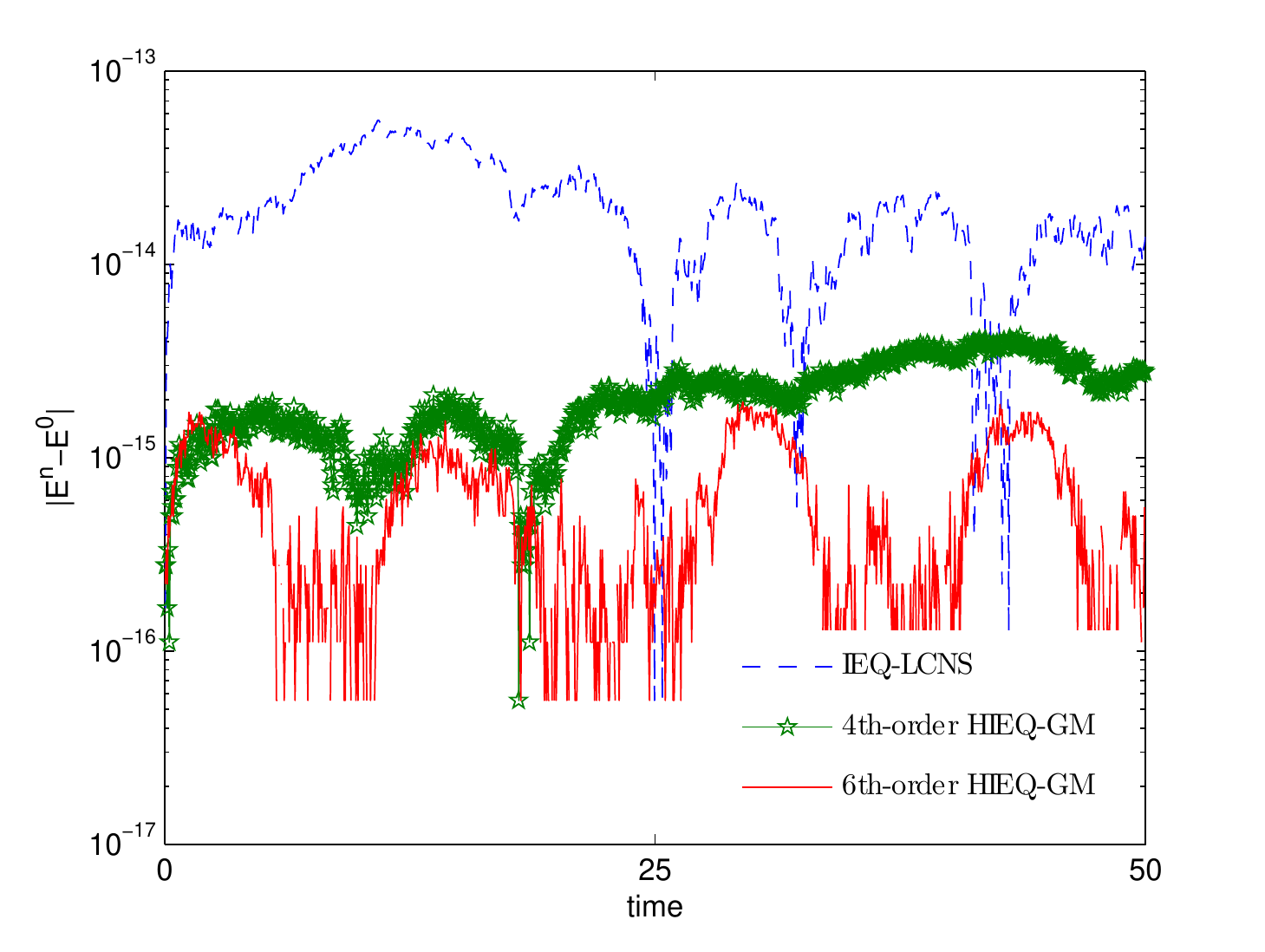}
\caption*{(d) Quadratic energy \eqref{ch-eq:4.7}}
\end{minipage}
\caption{The errors in invariants with $h=\frac{L}{128}$ and $\tau=0.0001$ over the time interval $t\in[0,50]$.}\label{Fig-ch:2}
\end{figure}

We consider the periodic peaked traveling wave with initial condition \cite{XS08}
 \begin{align*}
u_0(x)=\left\lbrace
  \begin{aligned}
  &\frac{c}{\cosh(L/2)}\cosh(x-x_0),\ |x-x_0|\le L/2,\\
  &\frac{c}{\cosh(L/2)}\cosh(L-(x-x_0)),\ |x-x_0|> L/2,
  \end{aligned}\right.\ \ 
  \end{align*}
where $c$ is the wave speed, $L$ is the period, and $x_0$ is the position of the trough. In the numerical experiment, 
the parameters are chosen as $c=1$, $L=1$, and $x_0=0$ and the periodic boundary condition is considered. The errors of invariants are plotted in Fig. \ref{Fig-ch:2}. In Fig. \ref{Fig-ch:2} (a), we can see that IEQ-LCNS and MSFPS can only preserve the Hamiltonian energy approximately and the error provided by IEQ-LCNS is largest. In theory, the proposed high-order schemes cannot exactly preserve the discrete Hamiltonian energy, however, {from Fig. \ref{Fig-ch:2} (a)}, we can observe that the resulting errors provided by 4th-order HIEQ-GM and 6th-order HIEQ-GM, respectively, can be preserved up to the machine accuracy and are much smaller  than the one provided by EPFPS. {In Figs. \ref{Fig-ch:2} (b)-(c) }, it is clear to see that the errors of the momentum are bounded and all of the schemes can exactly preserve the mass conservation law. { Fig. \ref{Fig-ch:2} (d)} show that the proposed schemes can exactly preserve the discrete quadratic energy \eqref{ch-eq:4.7}, which conforms the theoretical analysis.

\begin{figure}[H]
\centering\begin{minipage}[t]{60mm}
\includegraphics[width=65mm]{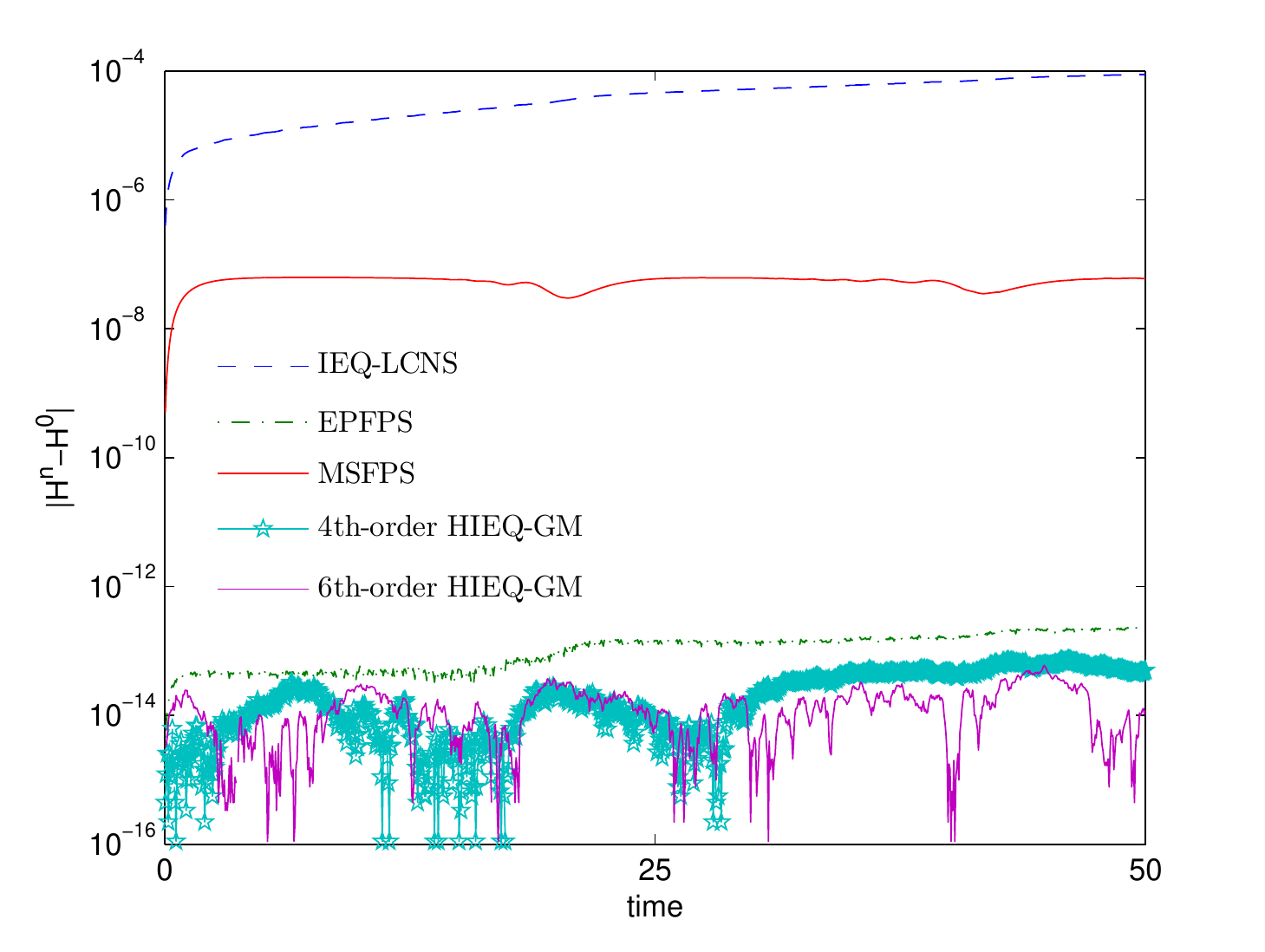}
\caption*{(a) Hamiltonian energy}
\end{minipage}\ \
\begin{minipage}[t]{60mm}
\includegraphics[width=65mm]{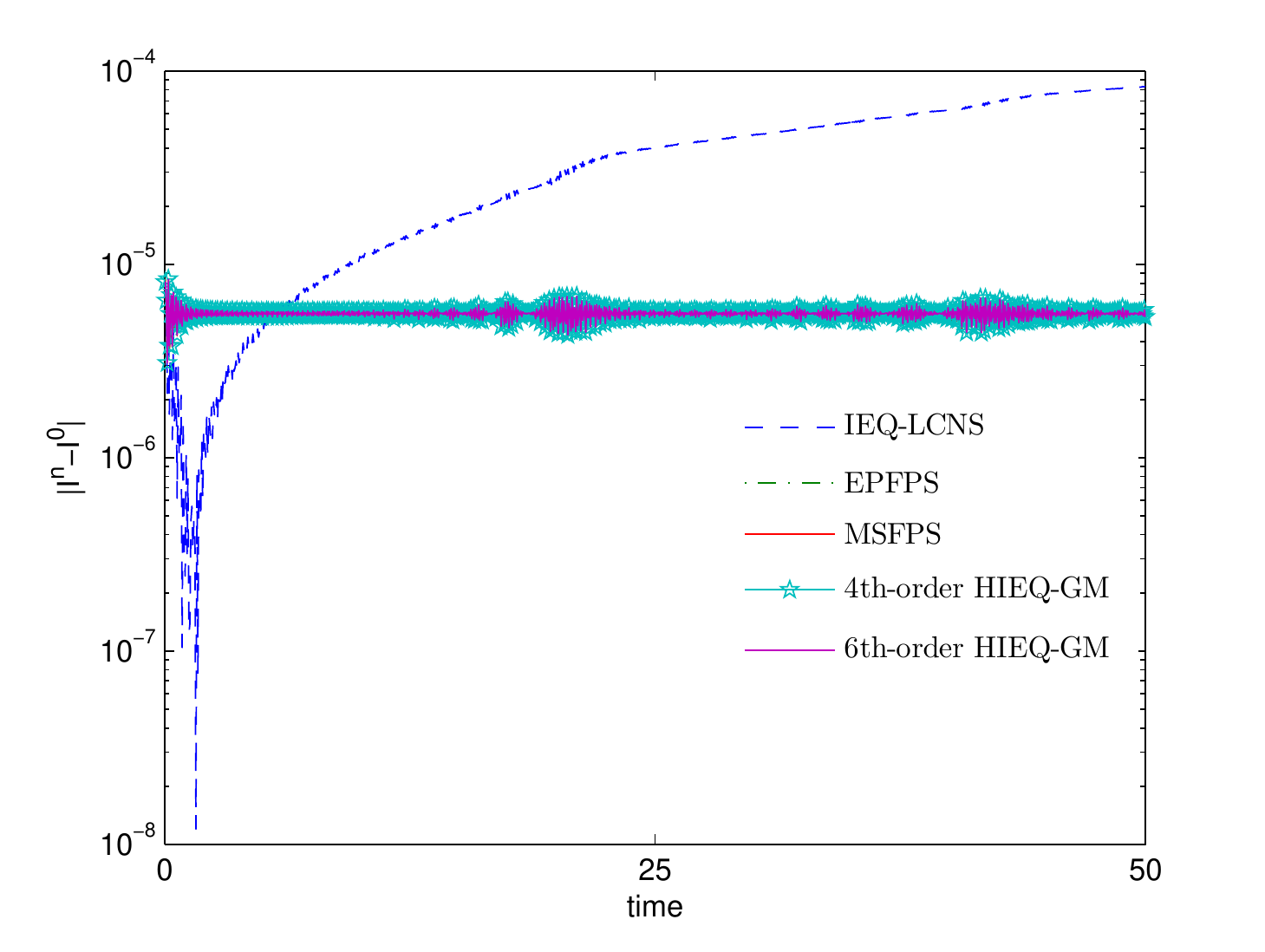}
\caption*{(b) Momentum}
\end{minipage}
\end{figure}

\subsection{Three-peakon interaction}

\begin{figure}[H]
\centering\begin{minipage}[t]{50mm}
\includegraphics[width=55mm]{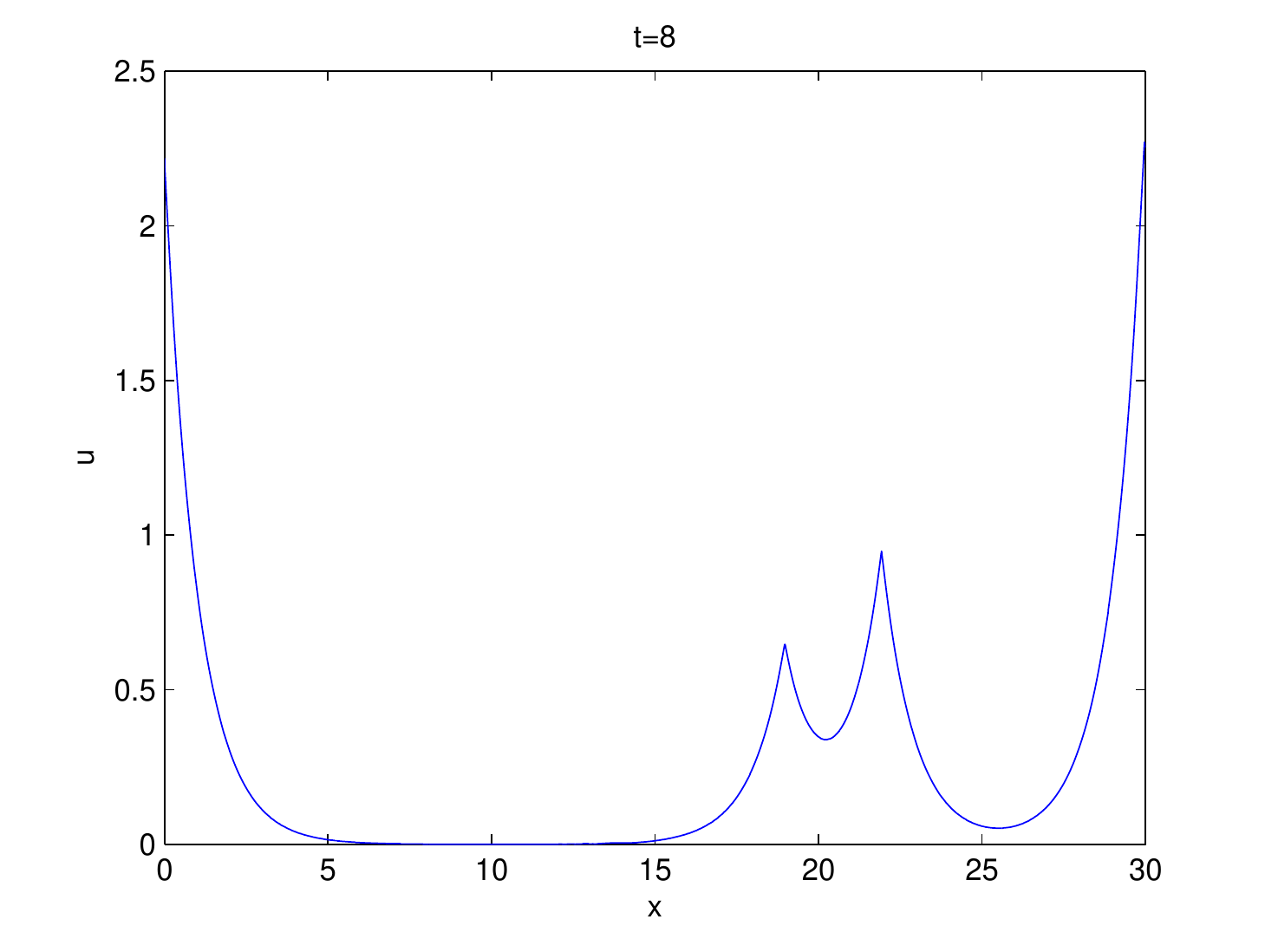}
\end{minipage}\ \
\begin{minipage}[t]{50mm}
\includegraphics[width=55mm]{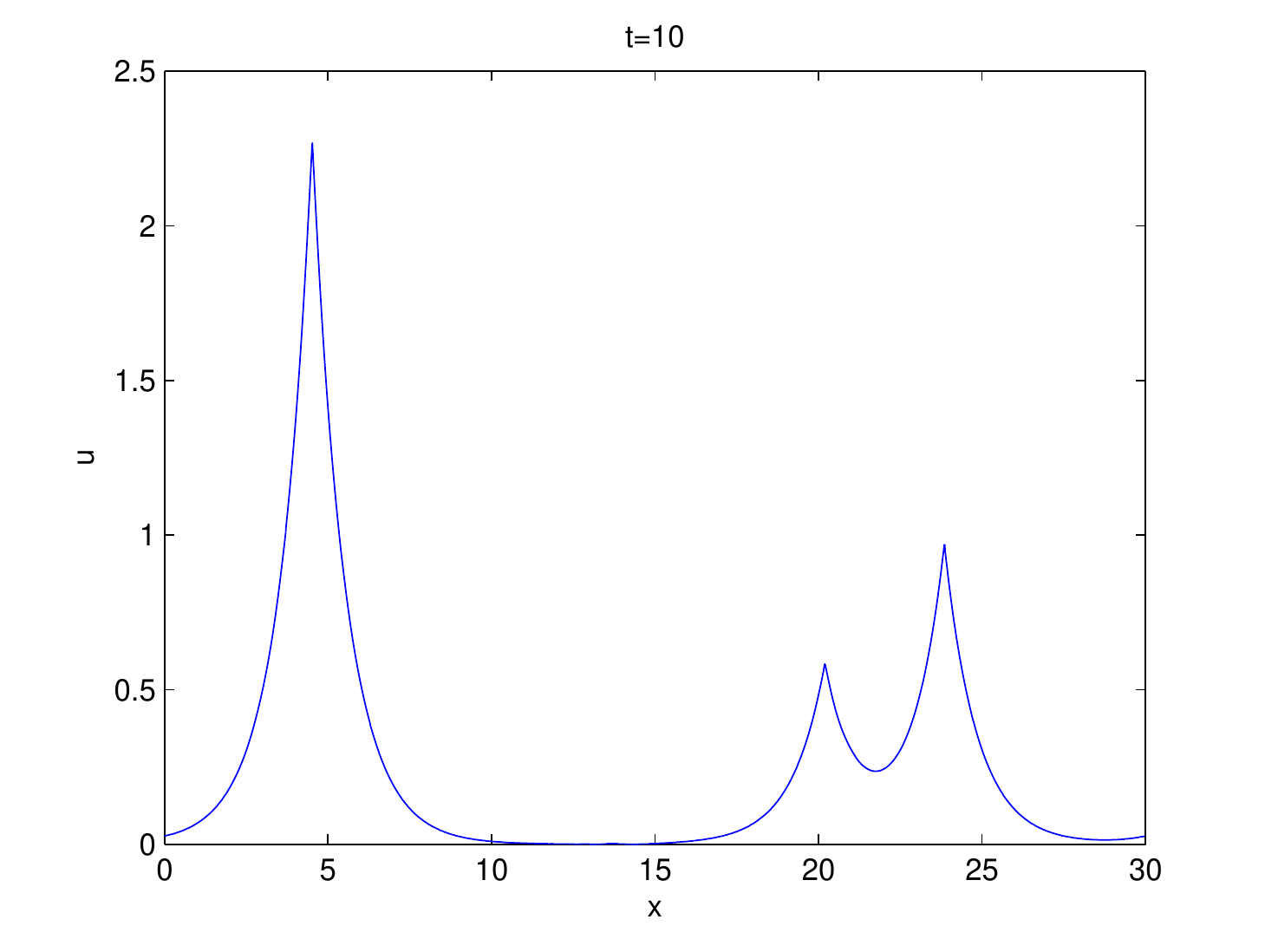}
\end{minipage}
\caption{The three-peakon interaction of the CH equation \eqref{ch-eq:1.1} provided by 4th-order HIEQ-GM with $h=\frac{L}{2048}$ and $\tau=0.0001$ at $t=0,1,2,3,4,6,8$ and $10$, respectively.}\label{Fig-ch:3}
\end{figure}

\begin{figure}[H]
\centering\begin{minipage}[t]{60mm}
\includegraphics[width=65mm]{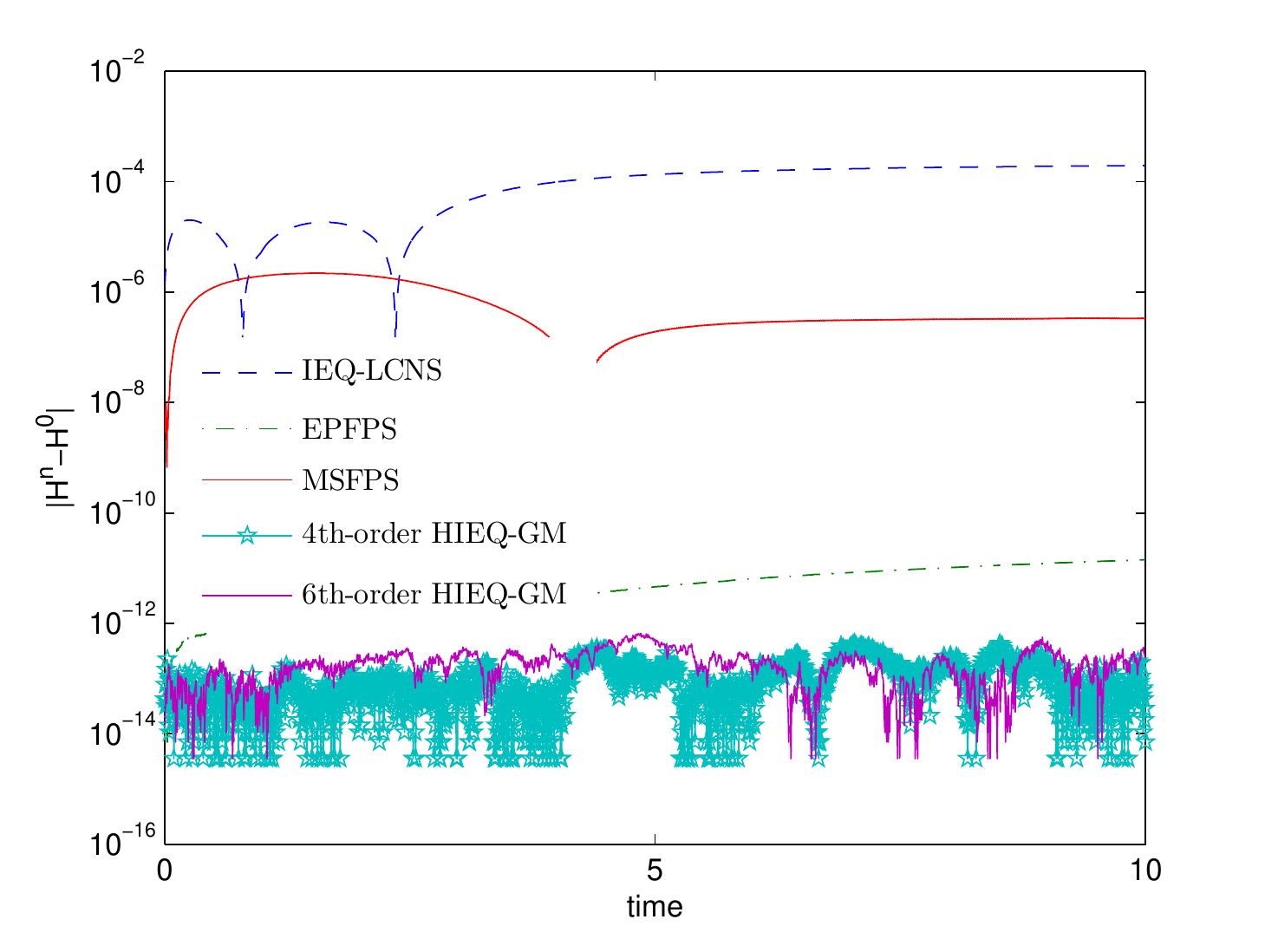}
\caption*{(a) Hamiltonian energy}
\end{minipage}\ \
\begin{minipage}[t]{60mm}
\includegraphics[width=65mm]{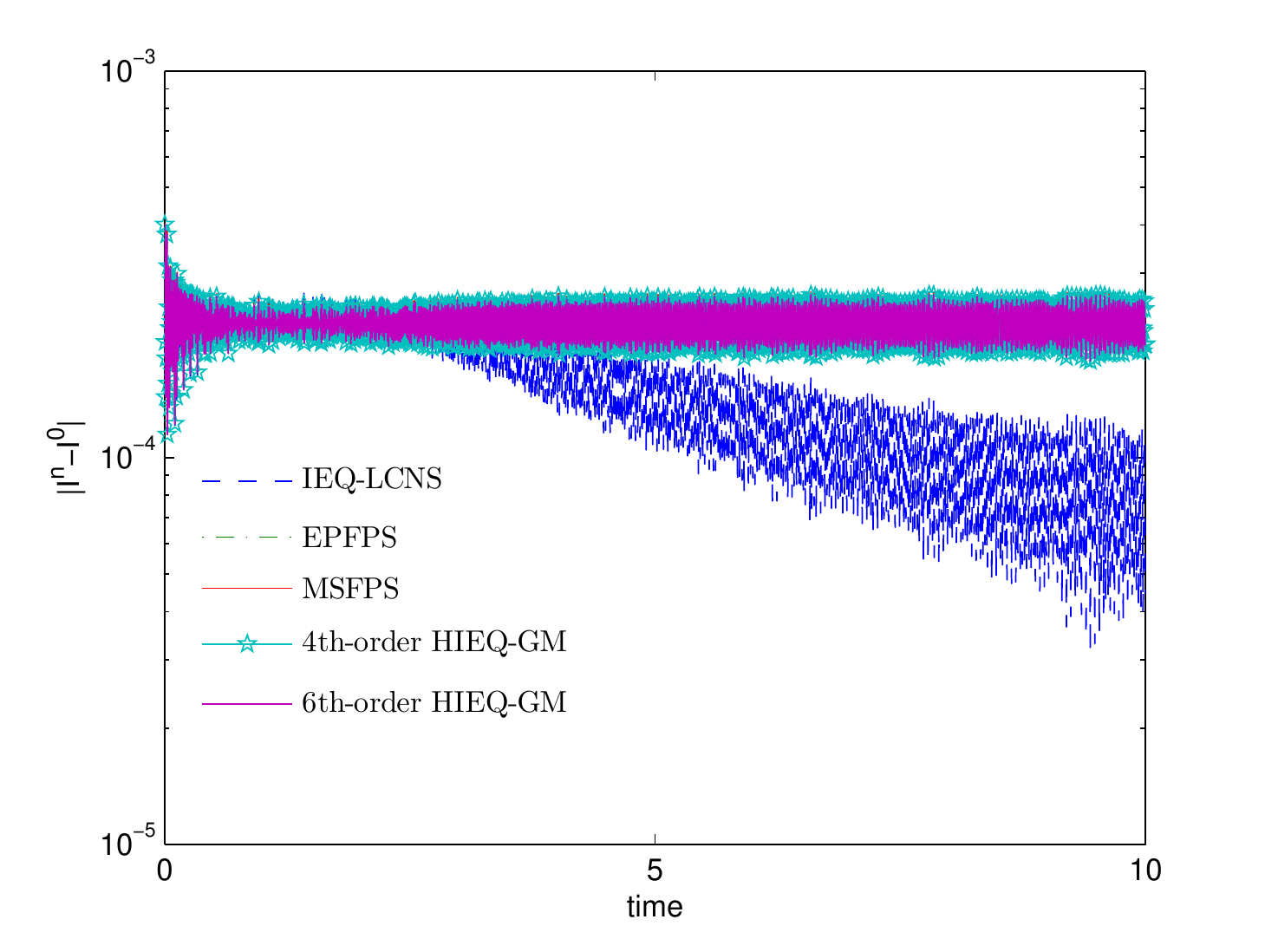}
\caption*{(b) Momentum}
\end{minipage}
\centering\begin{minipage}[t]{60mm}
\includegraphics[width=65mm]{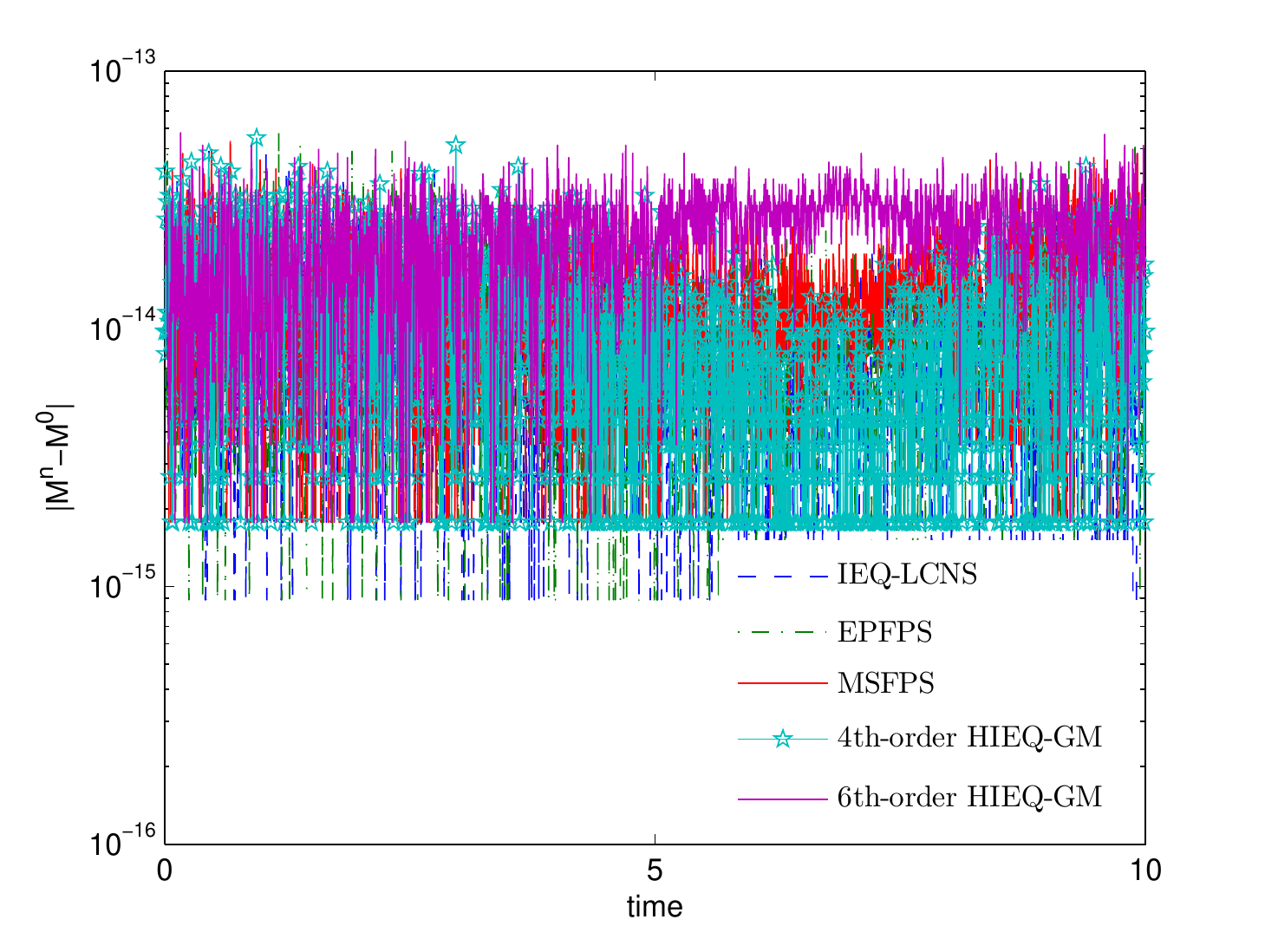}
\caption*{(c) Mass}
\end{minipage}\ \
\begin{minipage}[t]{60mm}
\includegraphics[width=65mm]{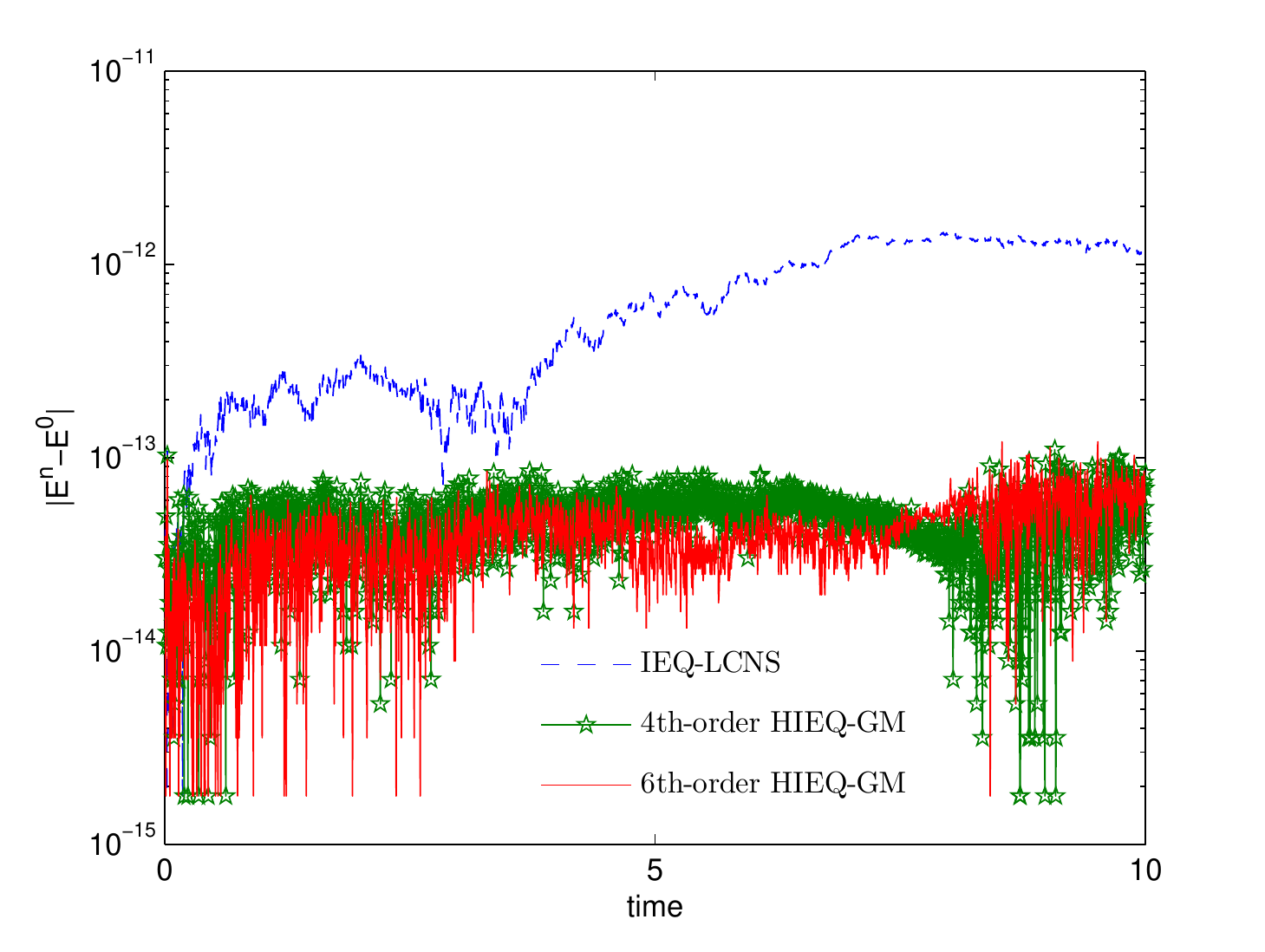}
\caption*{(d) Quadratic energy \eqref{ch-eq:4.7}}
\end{minipage}
\caption{The errors in invariants with $h=\frac{L}{2048}$ and $\tau=0.0001$ over the time interval $t\in[0,10]$.}\label{Fig-ch:4}
\end{figure}

In this example, we consider the three-peakon
interaction of the CH equation with initial condition \cite{XS08}
\begin{align*}
u_0(x)=\phi_1(x)+\phi_2(x)+\phi_3(x),
\end{align*}
where
\begin{align*}
\phi_i(x)=\left\lbrace
  \begin{aligned}
  &\frac{c_i}{\cosh(L/2)}\cosh(x-x_i),\ |x-x_i|\le L/2,\\
  &\frac{c_i}{\cosh(L/2)}\cosh(L-(x-x_i)),\ |x-x_i|> L/2,\
  \end{aligned} i=1,2,3.
  \right.\ \ 
  \end{align*}
The parameters are $c_1=2,c_2=1,c_3=0.8,x_1=-5,x_2=-3,x_3=-1$ and $L=30$, and the computational domain is $\Omega=[0,L]$ with the periodic
boundary condition. In Fig. \ref{Fig-ch:3}, we display the interaction of three peakons by 4th-order HIEQ-GM at $t=0,1,2,3,4,6,8$ and $10$, respectively. We can
see clearly that the moving peak interaction is resolved very well.
 The three-peakon interaction obtained by other schemes are not presented {since they are close to Fig. \ref{Fig-ch:3}}. The errors of invariants are plotted in Fig. \ref{Fig-ch:4}, which shows that all of the schemes can exactly preserve the mass conservation law and the momentum errors provided by the schemes are bounded. The Hamiltonian energy errors provided by the high-order schemes are smallest and the newly proposed schemes can exactly preserve the discrete quadratic energy \eqref{ch-eq:4.7}.

\begin{figure}[H]
\centering\begin{minipage}[t]{50mm}
\includegraphics[width=55mm]{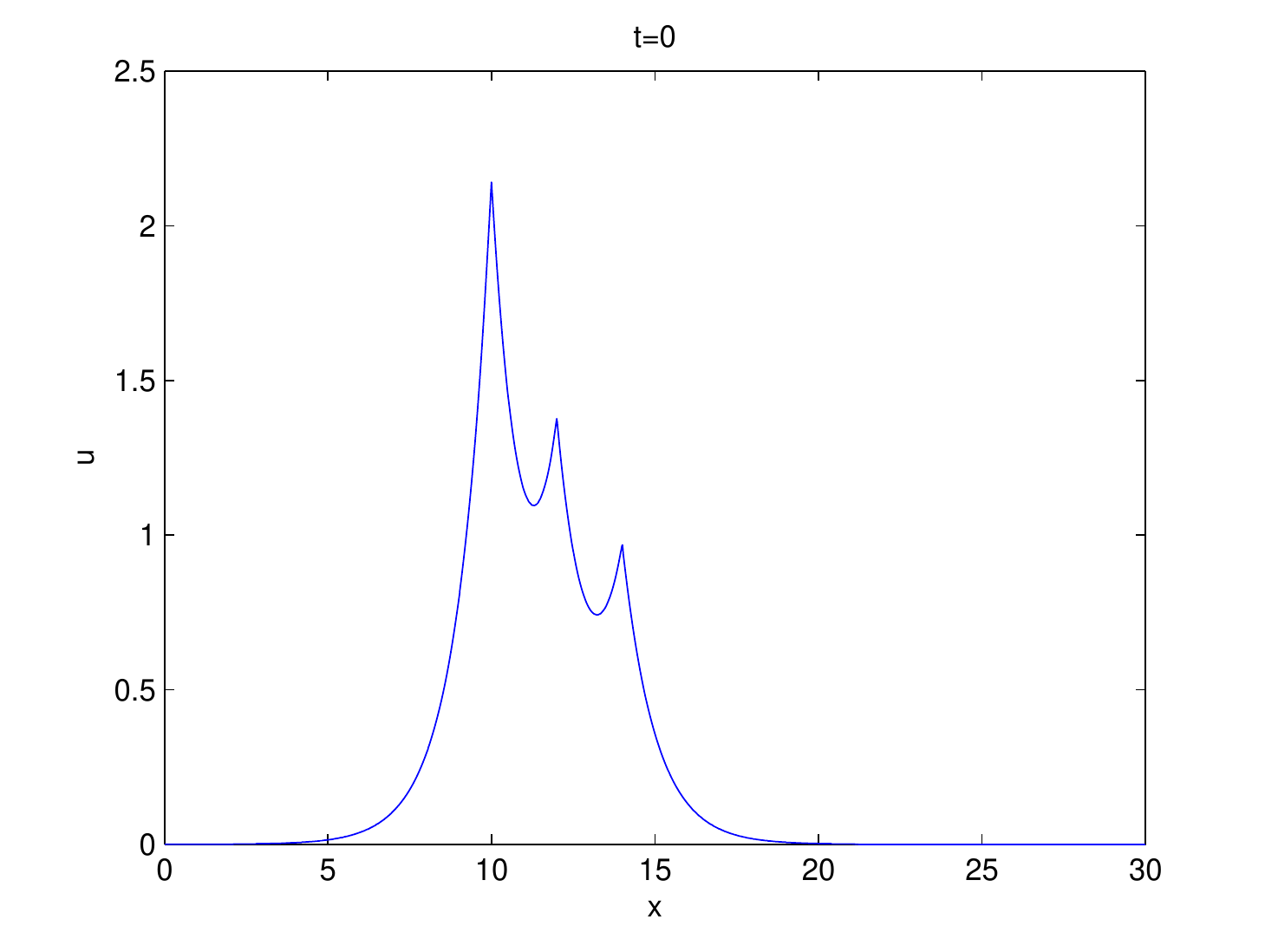}
\end{minipage}\ \
\begin{minipage}[t]{50mm}
\includegraphics[width=55mm]{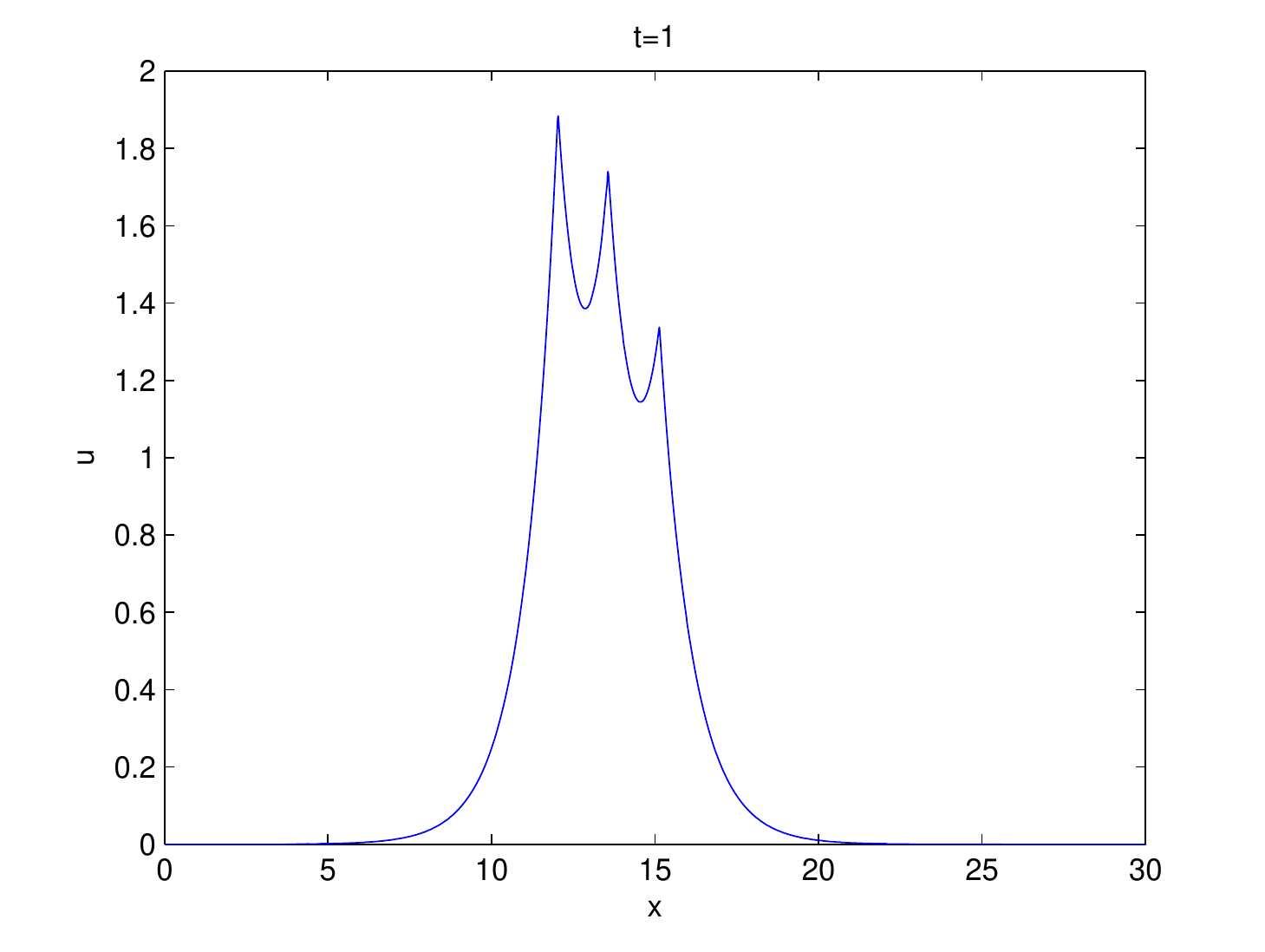}
\end{minipage}
\centering\begin{minipage}[t]{50mm}
\includegraphics[width=55mm]{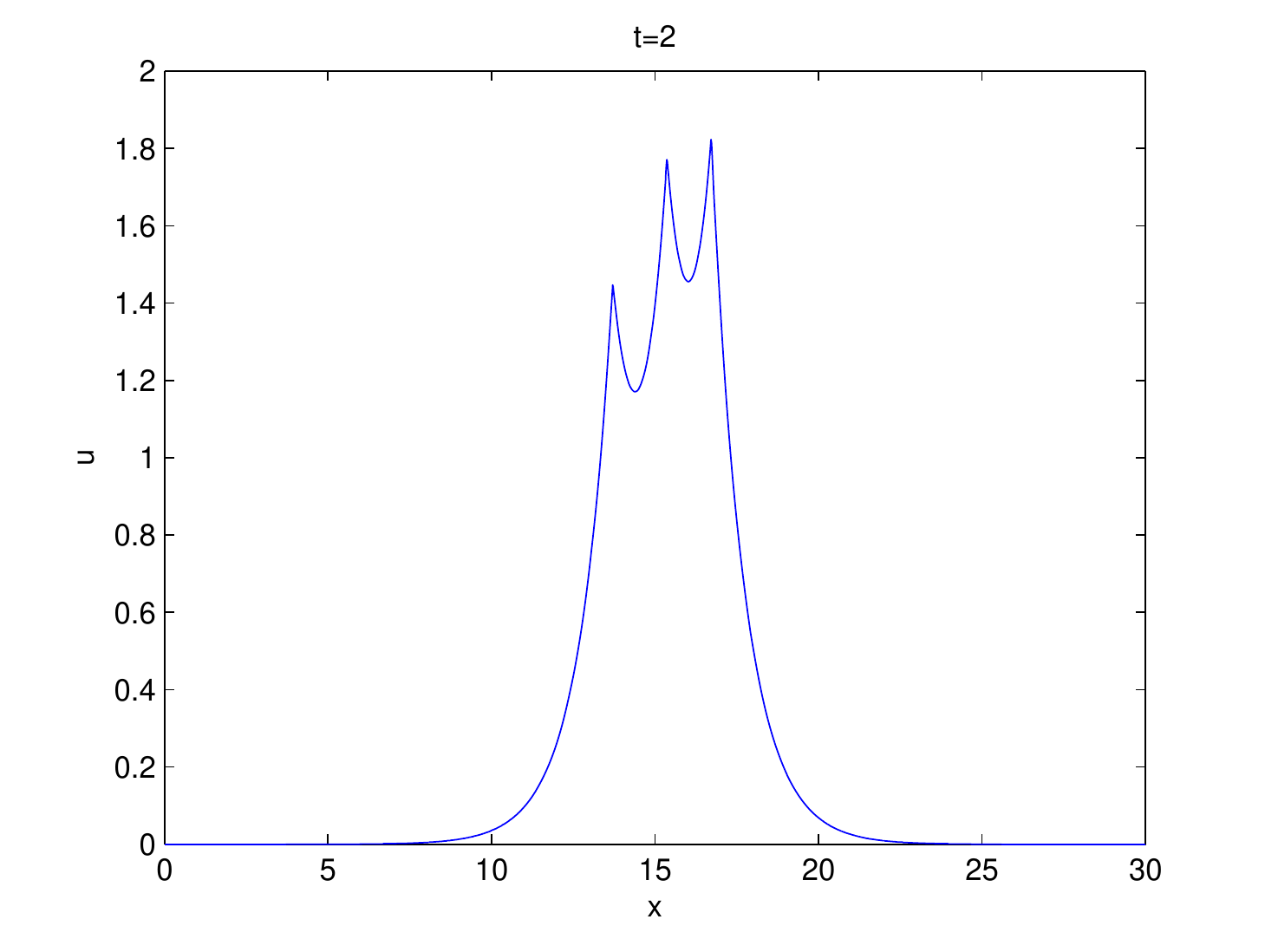}
\end{minipage}\ \
\begin{minipage}[t]{50mm}
\includegraphics[width=55mm]{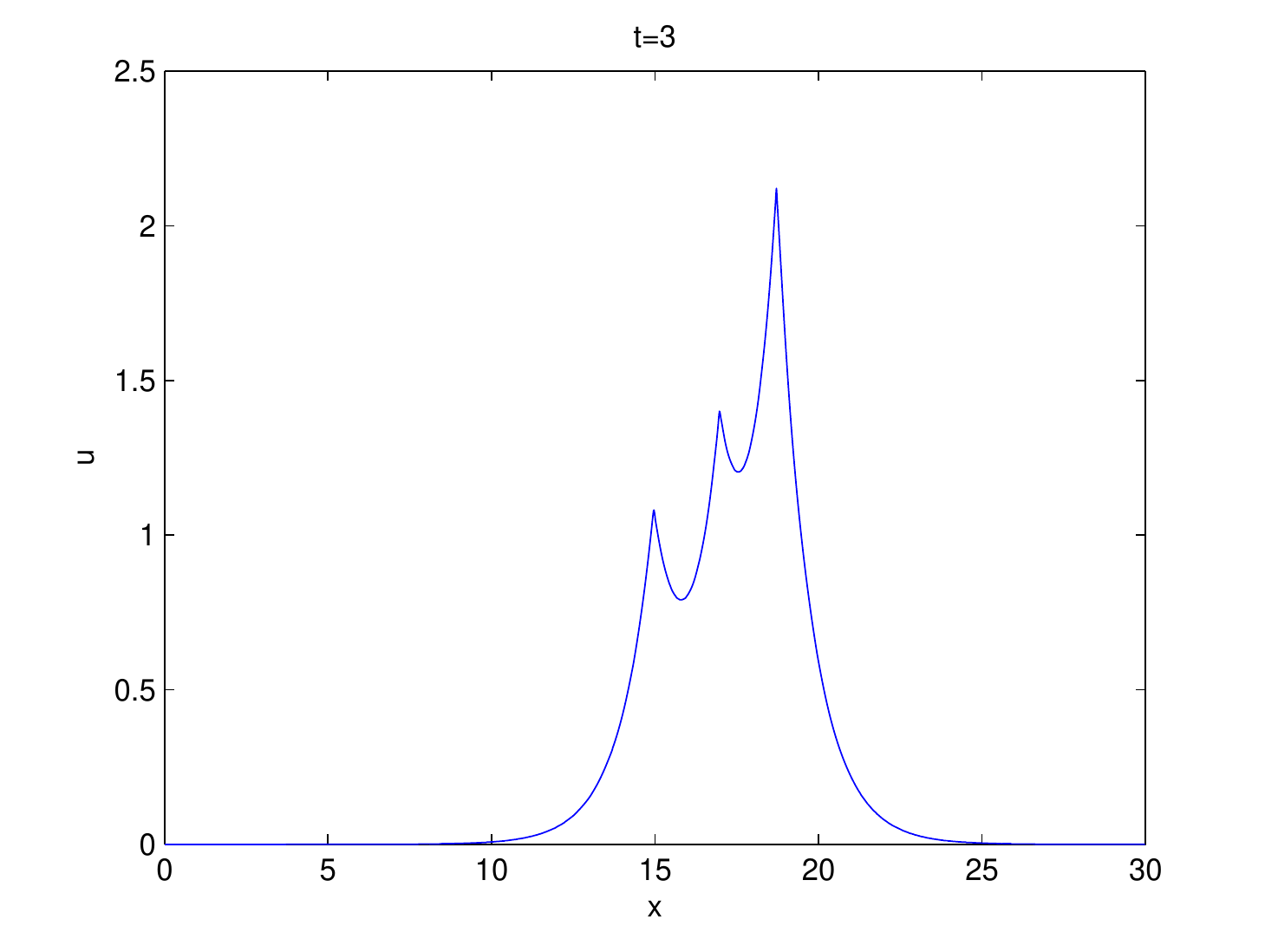}
\end{minipage}
\centering\begin{minipage}[t]{50mm}
\includegraphics[width=55mm]{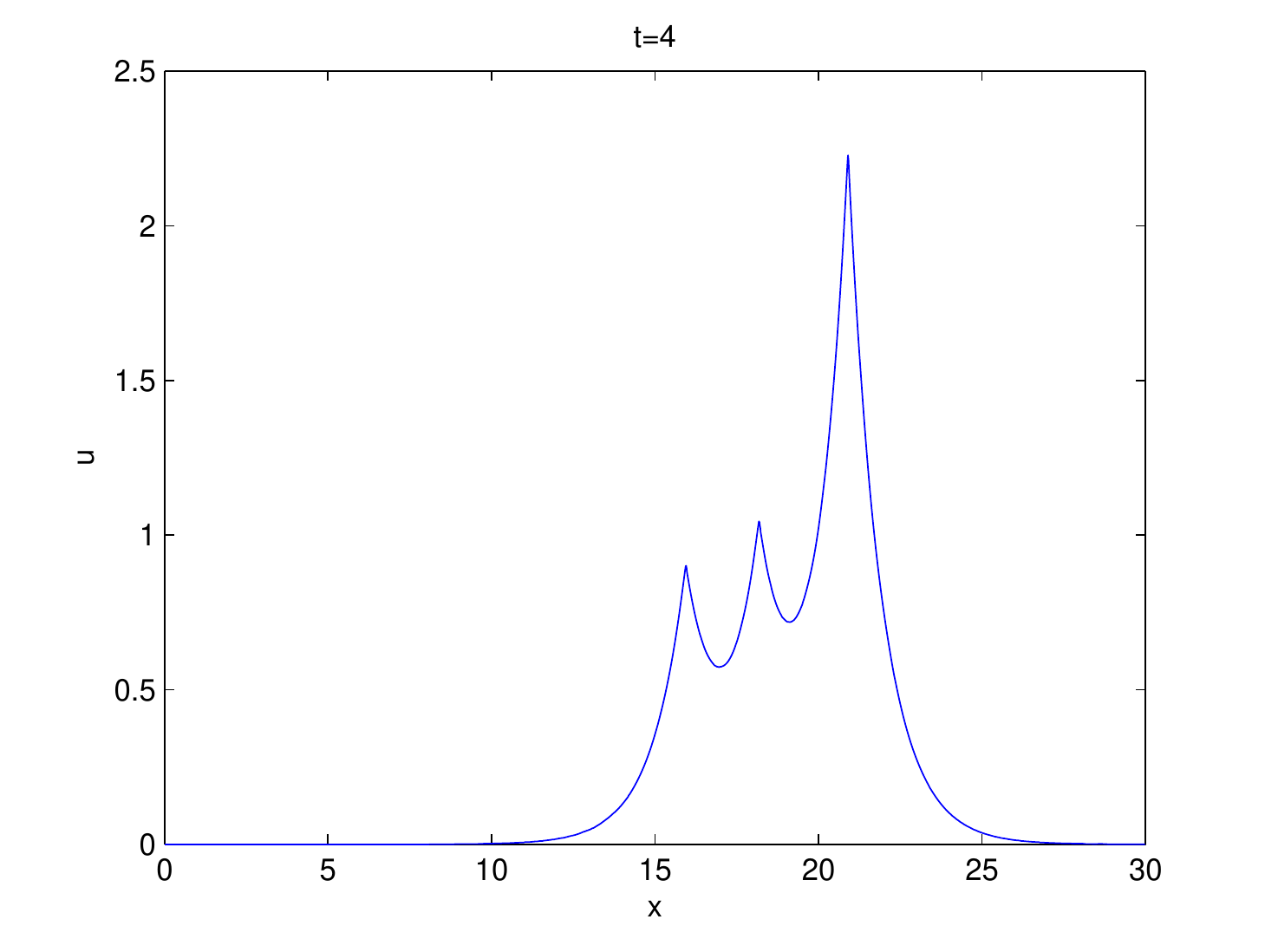}
\end{minipage}\ \
\begin{minipage}[t]{50mm}
\includegraphics[width=55mm]{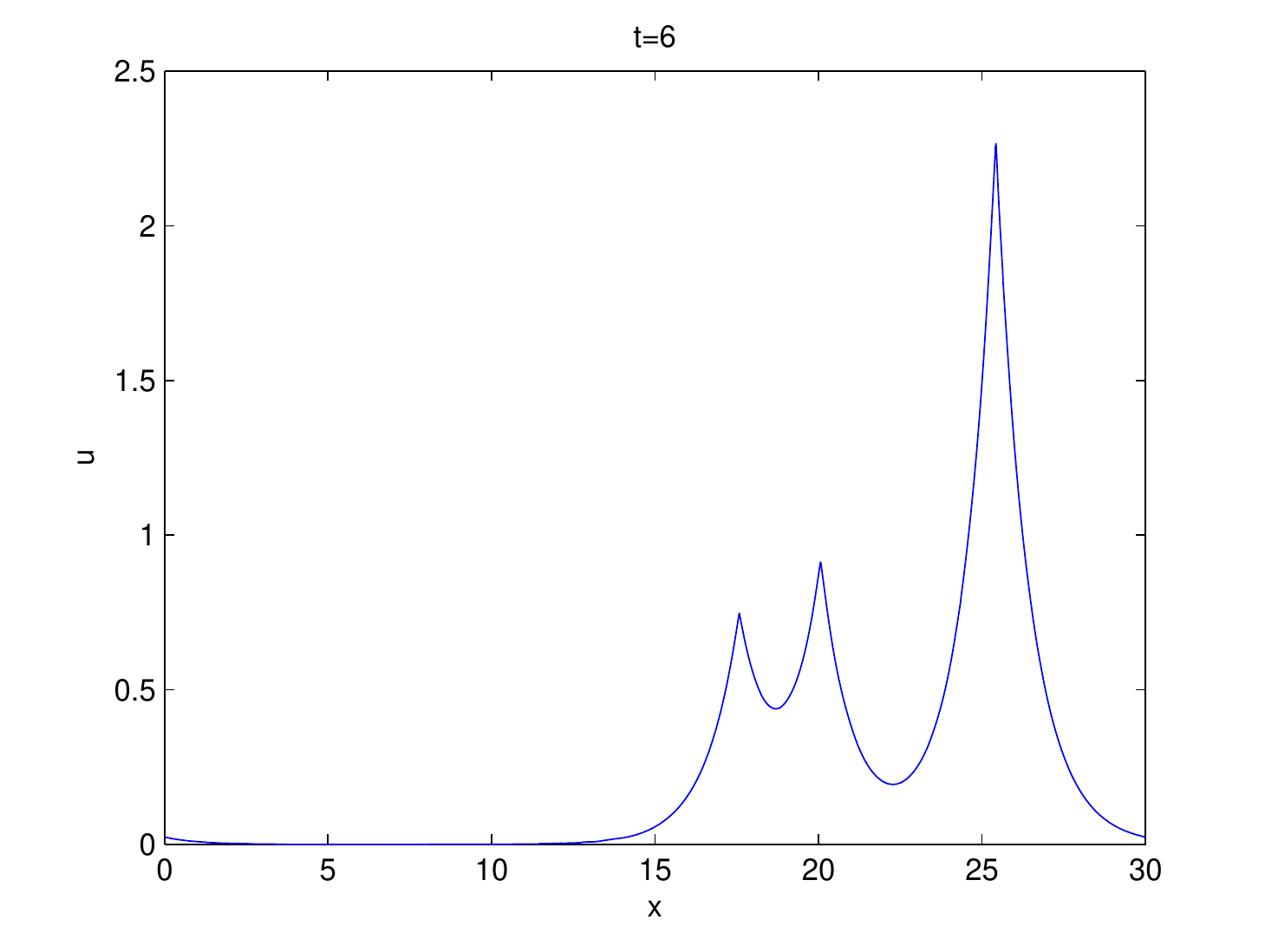}
\end{minipage}
\end{figure}

\section{Concluding remarks}\label{ch:Sec.6}
In this paper, we combine the idea of the IEQ approach with symplectic Runge-Kutta methods to
propose a new class of energy-preserving methods for the CH equation. The proposed schemes could reach arbitrarily high-order accuracy while exactly preserving
the discrete quadratic energy of the modified system. Numerical examples are addressed to illustrate the accuracy and
energy-preserving property of the proposed schemes. Compared with some existing low-order structure-preserving schemes, the proposed high-order
schemes show remarkable efficiency and the advantage in preserving the discrete Hamiltonian energy.

{We conclude this paper with some remarks. First, the presented strategy can be directly extended to propose high-order energy-preserving methods for the  Hamiltonian PDEs where the Hamiltonian functionals are quadratic. For a general case (including the non-canonical Hamiltonian PDE), we dealt with it, as follows: we first utilize the idea of the IEQ approach to transform Hamiltonian energy as a quadratic form and then, following the energy variational, the original system is reformulated into an equivalent system, which inherits such quadratic energy. Finally, the resulting system is solved by using in time a symplectic Runge-Kutta method. 
Second, compared with a existing energy-preserving method such as HBVMs, the proposed method can not preserve the
Hamiltonian energy of the original system. Thus, such trade-offs among methods should be more carefully investigated. Finally, we can also apply the idea of the scalar auxiliary variable (SAV) approach \cite{SXY17,SXY18} to reformulate the CH equation into a new equivalent system which inherits a quadratic energy. Thus, a possible future work is to develop high-order energy-preserving methods for Hamiltonian PDEs by combine the idea of the SAV approach with the symplectic Runge-Kutta method.}

\section*{Acknowledgments} The authors would like to express sincere gratitude to the referees for their insightful
comments and suggestions. Chaolong Jiang's work is partially supported by the National Natural Science Foundation of China (Grant No. 11901513), the Yunnan Provincial Department of Education Science Research Fund Project (Grant No. 2019J0956) and the Science and Technology Innovation Team on Applied Mathematics in Universities of Yunnan. Yushun Wang's work is partially supported by the National Natural Science Foundation of China (Grant No. 11771213) and the National Key Research and Development Project of China (Grant Nos. 2016YFC0600310, 2018YFC0603500, 2018YFC1504205). Yuezheng Gong's work is partially supported by the Natural Science Foundation of Jiangsu Province
(Grant No. BK20180413) and the National Natural Science Foundation of China (Grant No. 11801269).


\begin{thebibliography}{10}

\bibitem{BCMR12}
L.~Brugnano, M.~Calvo, J. I. Montijano, and L.~{R\'andez}.
\newblock Energy-preserving methods for {P}oisson systems.
\newblock {\em J. Comput. Appl. Math.}, 236:3890--3904, 2012.

\bibitem{BI16}
L.~Brugnano and F.~Iavernaro.
\newblock {\em Line Integral Methods for Conservative Problems}.
\newblock Chapman et Hall/CRC: Boca Raton, FL, USA, 2016.

\bibitem{BIT10}
L.~Brugnano, F.~Iavernaro, and D.~Trigiante.
\newblock Hamiltonian boundary value methods (energy preserving discrete line
  integral methods).
\newblock {\em J. Numer. Anal. Ind. Appl. Math.}, 5:17--37, 2010.

\bibitem{CHWG15}
J. Cai, J. Hong, Y. Wang, and Y. Gong.
\newblock Two energy-conserved splitting methods for three-dimensional
  time-domain {M}axwell's equations and the convergence analysis.
\newblock {\em SIAM. J. Numer. Anal.}, 53:1918--1940, 2015.

\bibitem{CSW16}
W. Cai, Y. Sun, and Y. Wang.
\newblock Geometric numerical integration for peakon b-family equations.
\newblock {\em Commun. Comput. Phys.}, 19:24--52, 2016.

\bibitem{CH93}
R.~Camassa and D.~Holm.
\newblock An integrable shallow water equation with peaked solitons.
\newblock {\em Phys. Rev. Lett.}, 71:1661--1664, 1993.

\bibitem{CHH94}
R.~Camassa, D.~Holm, and J.~Hyman.
\newblock A new integrable shallow water equation.
\newblock {\em Adv. Appl. Mech.}, 31:1--33, 1994.

\bibitem{CL08}
R.~Camassa and L.~Lee.
\newblock Complete integrable particle methods and the recurrence of initial
  states for a nonlinear shallow-water wave equation.
\newblock {\em J. Comput. Phys.}, 227:7206--7221, 2008.

\bibitem{CQ01}
J. Chen and M. Qin.
\newblock Multi-symplectic {F}ourier pseudospectral method for the nonlinear
  {S}chr\"{o}dinger equation.
\newblock {\em Electr. Trans. Numer. Anal.}, 12:193--204, 2001.

\bibitem{CKR08a}
G.~Coclite, K.~Karlsen, and N.~Risebro.
\newblock A convergent finite difference scheme for the {C}amassa-{H}olm
  equation with general {$H^1$} initial data.
\newblock {\em SIAM J. Numer. Anal.}, 46:1554--1579, 2008.

\bibitem{CH11}
D.~Cohen and E.~Hairer.
\newblock Linear energy-preserving integrators for {P}oisson systems.
\newblock {\em BIT}, 51:91--101, 2011.

\bibitem{COR08}
D.~Cohen, B.~Owren, and X.~Raynaud.
\newblock Multi-symplectic integration of the {C}amassa-{H}olm equation.
\newblock {\em J. Comput. Phys.}, 227:5492--5512, 2008.

\bibitem{CR11}
D.~Cohen and X.~Raynaud.
\newblock Geometric finite difference schemes for the generalized
  hyperelastic-rod wave equation.
\newblock {\em J. Comput. Appl. Math.}, 235:1925--1940, 2011.

\bibitem{CE98}
A.~Constantin and J.~Escher.
\newblock Global existence and blow-up for a shallow water equation.
\newblock {\em Ann. Scuola Norm. Sup. Pisa Cl. Sci.}, 26:303--328, 1998.

\bibitem{Cooperima87}
{G.~J. Cooper.
\newblock Stability of Runge-Kutta methods for trajectory problems.
\newblock {\em IMA J. Numer. Anal.}, 7:1--13, 1987.}

\bibitem{ELS19}
S.~Eidnes, L.~Li, and S.~Sato.
\newblock Linearly implicit structure-preserving schemes for {H}amiltonian
  systems.
\newblock {\em arXiv preprint arXiv:1901.03573}, 2019.

\bibitem{FL09}
B. Feng and Y.~Liu.
\newblock An operator splitting method for the {D}egasperis-{P}rocesi equation.
\newblock {\em J. Comput. Phys.}, 228:7805--7820, 2009.

\bibitem{FMO10}
B. Feng, K.~Maruno, and Y.~Ohta.
\newblock A self-adaptive moving mesh method for the {C}amassa-{H}olm equation.
\newblock {\em J. Comput. Appl. Math.}, 235:229--243, 2010.

\bibitem{GW16b}
Y. Gong and Y. Wang.
\newblock An energy-preserving wavelet collocation method for general
  multi-symplectic formulations of {H}amiltonian {PDE}s.
\newblock {\em Commun. Comput. Phys.}, 20:1313--1339, 2016.

\bibitem{GWW18}
Y. Gong, Y. Wang, and Q.~Wang.
\newblock Linear-implicit conservative schemes based on energy quadratization
  for {H}amiltonian {PDE}s.
\newblock {\em Preprint.}

\bibitem{GZYW18}
Y. Gong, J.~Zhao, X. Yang, and Q.~Wang.
\newblock Fully discrete second-order linear schemes for hydrodynamic phase
  field models of binary viscous fluid flows with variable densities.
\newblock {\em SIAM J. Sci. Comput.}, 40:B138--B167, 2018.

\bibitem{H10}
E.~Hairer.
\newblock Energy-preserving variant of collocation methods.
\newblock {\em J. Numer. Anal. Ind. Appl. Math.}, 5:73--84, 2010.

\bibitem{ELW06}
E.~Hairer, C.~Lubich, and G.~Wanner.
\newblock {\em Geometric Numerical Integration: Structure-Preserving Algorithms
  for Ordinary Differential Equations}.
\newblock Springer-Verlag, Berlin, 2nd edition, 2006.

\bibitem{HR06}
H.~Holden and X.~Raynaud.
\newblock Convergence of a finite difference scheme for the {C}amassa-{H}olm
  equation.
\newblock {\em SIAM J. Numer. Anal.}, 44:1655--1680, 2006.

\bibitem{HGL19}
Q.~Hong, Y. Gong, and Z. Lv.
\newblock Linear and {H}amiltonian-conserving {F}ourier pseudo-spectral schemes
  for the {C}amassa-{H}olm equation.
\newblock {\em Appl. Math. Comput.}, 346:86--95, 2019.

\bibitem{JCW18}
C. Jiang, W. Cai, and Y. Wang.
\newblock A linearly implicit and local energy-preserving scheme for the sine-{G}ordon equation based on the invariant energy quadratization approach.
\newblock {\em J. Sci. Comput.}, 80:1629-1655, 2019.

\bibitem{JCWL17}
C. Jiang, W. Cai, Y. Wang, and H. Li.
\newblock A sixth order energy-conserved method for three-dimensional
  time-domain {M}axwell's equations.
\newblock {\em arXiv preprint}, \:arXiv:1705.08125, 2017.

\bibitem{KL05}
H.~Kalisch and J.~Lenells.
\newblock Numerical study of traveling-wave solutions for the {C}amassa- {H}olm
  equation.
\newblock {\em Chaos Solitons Fractals}, 25:287--298, 2005.

\bibitem{LO00}
A.~Li and P.~Olver.
\newblock Well-posedness and blow-up solutions for an integrable nonlinearly
  dispersive model wave equation.
\newblock {\em J. Differ. Equ.}, 162:27--63, 2000.

\bibitem{LWQ14}
H. Li, Y. Wang, and M. Qin.
\newblock A sixth order averaged vector field method.
\newblock {\em J. Comput. Math.}, 34:479--498, 2016.

\bibitem{LW16}
Y. Li and X. Wu.
\newblock Functionally fitted energy-preserving methods for solving oscillatory
  nonlinear {H}amiltonian systems.
\newblock {\em SIAM J. Numer. Anal.}, 54:2036--2059, 2016.

\bibitem{Matsuo10}
T.~Matsuo.
\newblock A {H}amiltonian-conserving {G}alerkin scheme for the {C}amassa-{H}olm
  equation.
\newblock {\em J. Comput. Appl. Math.}, 234:1258--1266, 2010.

\bibitem{MY09}
T.~Matsuo and H.~Yamaguchi.
\newblock An energy-conserving {G}alerkin scheme for a class of nonlinear
  dispersive equations.
\newblock {\em J. Comput. Phys.}, 228:4346--4358, 2009.

\bibitem{Miyatake14}
Y.~Miyatake.
\newblock An energy-preserving exponentially-fitted continuous stage
  {R}unge-{K}utta method for {H}amiltonian systems.
\newblock {\em BIT}, 54:777--799, 2014.

\bibitem{QM08}
G. R. W. Quispel and D. I. McLaren.
\newblock A new class of energy-preserving numerical integration methods.
\newblock {\em J. Phys. A: Math. Theor.}, 41:045206, 2008.

\bibitem{Sanz-Sernabit88}
{J.~M. Sanz-Serna.
\newblock Runge-Kutta schemes for Hamiltonian systems.
\newblock {\em BIT}, 28:877--883, 1988.}

\bibitem{SCbook94}
{J.~M. Sanz-Serna and M.~P. Calvo.
\newblock {\em Numerical Hamiltonian Problems}.
\newblock Chapman \& Hall, London, 1994.}

\bibitem{ST06}
J.~Shen and T.~Tang.
\newblock {\em Spectral and High-Order Methods with Applications}.
\newblock Science Press, Beijing, 2006.

\bibitem{SXY17}
J.~Shen, J.~Xu, and J.~Yang.
\newblock A new class of efficient and robust energy stable schemes for
  gradient flows.
\newblock {\em SIAM Rev.}, 61:474--506, 2019.


\bibitem{SXY18}
J.~Shen, J.~Xu, and J.~Yang.
\newblock The scalar auxiliary variable {(SAV)} approach for gradient.
\newblock {\em J. Comput. Phys.}, 353:407--416, 2018.

\bibitem{TS12}
W. Tang and Y. Sun.
\newblock Time finite element methods: a unified framework for numerical
  discretizations of {ODE}s.
\newblock {\em Appl. Math. Comput.}, 219:2158--2179, 2012.

\bibitem{WW18}
B.~Wang and X. Wu.
\newblock Functionally-fitted energy-preserving integrators for {P}oisson
  systems.
\newblock {\em J. Comput. Phys.}, 364:137--152, 2018.

\bibitem{XS08}
Y.~Xu and C.-W Shu.
\newblock A local discontinuous {G}alerkin method for the {C}amassa-{H}olm
  equation,.
\newblock {\em SIAM J. Numer. Anal.}, 46:1998--2021, 2008.

\bibitem{YZW17}
X. Yang, J.~Zhao, and Q.~Wang.
\newblock Numerical approximations for the molecular beam epitaxial growth
  model based on the invariant energy quadratization method.
\newblock {\em J. Comput. Phys.}, 333:104--127, 2017.

\bibitem{YZWS17}
X. Yang, J.~Zhao, Q.~Wang, and J.~Shen.
\newblock Numerical approximations for a three components {C}ahn-{H}illiard
  phase-field model based on the invariant energy quadratization method.
\newblock {\em Math. Models Methods Appl. Sci.}, 27:1993--2030, 2017.

\bibitem{ZYGW17}
J.~Zhao, X. Yang, Y. Gong, and Q.~Wang.
\newblock A novel linear second order unconditionally energy stable scheme for
  a hydrodynamic-tensor model of liquid crystals.
\newblock {\em Comput. Methods Appl. Mech. Engrg.}, 318:803--825, 2017.

\bibitem{ZST11}
H. Zhu, S. Song, and Y. Tang.
\newblock Multi-symplectic wavelet collocation method for the {S}chr\"{o}dinger
  equation and the {C}amassa-{H}olm equation.
\newblock {\em Comput. Phys. Commun.}, 182:616--627, 2011.

\end{thebibliography}
\end{document}